\documentclass{article}
\usepackage{arxiv}
\usepackage[utf8]{inputenc}
\usepackage[T1]{fontenc}
\usepackage{amsmath,amssymb,amsthm}
\usepackage{graphicx,booktabs,subcaption}
\usepackage{placeins}
\usepackage{microtype}
\usepackage[numbers,sort&compress]{natbib}
\usepackage{hyperref}

\usepackage{algorithm}
\usepackage{algpseudocode}
\newcommand{\R}{\mathbb{R}}
\newcommand{\T}{\mathsf{T}}
\newcommand{\avg}[1]{\{\!\{#1\}\!\}}
\newcommand{\jump}[1]{[\![#1]\!]}
\title{A Novel Matrix-Free Approximate Elemental Inverse for Efficient Preconditioning of Discontinuous Galerkin Spectral Elements}
\author{
  Boyang Xia \\
  Department of Engineering \\
  King's College London \\
  London, United Kingdom \\
  \texttt{boyang.xia@kcl.ac.uk}
  \And
  David Moxey \thanks{Corresponding author.} \\
  Department of Engineering \\
  King's College London \\
  London, United Kingdom \\
  \texttt{david.moxey@kcl.ac.uk}
}
\date{13 September 2026}
\renewcommand{\shorttitle}{Matrix-free approximate elemental inverse}
\hypersetup{hidelinks,pdftitle={A Novel Matrix-Free Approximate Elemental Inverse for Efficient Preconditioning of Discontinuous Galerkin Spectral Elements},pdfsubject={Numerical methods for elliptic PDEs},pdfkeywords={discontinuous Galerkin, spectral element, approximate inverse, multigrid, matrix-free}}
\begin{document}
\raggedbottom
\maketitle
\begin{abstract}
Efficient matrix-free preconditioning remains a challenge for discontinuous
Galerkin spectral element discretisations of elliptic problems, particularly
for general element shapes with collapsed-coordinate bases, such as prisms,
tetrahedra and pyramids. Multigrid provides a
promising framework, but its efficiency depends strongly on the effectiveness
and computational cost of the smoother. For modal bases, point Jacobi can
provide insufficient smoothing, whereas applying block Jacobi through
stored dense inverses incurs substantial storage and computational costs.
We propose a matrix-free approximate elemental inverse, an approximation to the
inverse of each elemental diagonal block, constructed from
inverse modal-to-nodal transforms and an element-dependent sparse matrix.
Its application exploits generalised sum-factorisation and partial
sum-factorisation in collapsed coordinates, avoiding the explicit storage
and multiplication of large dense matrices. The method is evaluated as
both a standalone preconditioner and a multigrid smoother on hexahedral,
prismatic, tetrahedral and pyramidal elements. Although the proposed approximation
generally increases iteration counts compared to exact block Jacobi,
the cost per iteration is far lower. This trade-off yields
substantial reductions in total solution time for hexahedral and prismatic
elements, while achieving comparable performance for tetrahedral and
pyramidal elements. This study can also serve as a practical guide for tuning 
the convergence of the discontinuous Galerkin interior-penalty Poisson solver.
\end{abstract}


\section{Introduction}
\label{sec:introduction}
Matrix-free implementations and sum-factorisation are established approaches
to efficient high-order continuous or discontinuous Galerkin (DG) spectral element discretisations~\citep{kronbichler2019,kempf2020AutoCodeGen,fischer2022nekrs,jansson2024neko,andrej2024mfem}.
By evaluating all the integral and flux terms on the fly, matrix-free implementations avoid storing and accessing the assembled elemental or global sparse system matrix, which significantly reduces memory traffic and increases the arithmetic intensity.
Sum-factorisation further reduces the arithmetic work by expressing multidimensional elemental
operators as a series of one-dimensional tensor contractions. This can reduce both data access and floating-point operations. As modern computing hardware is well suited to algorithms with high arithmetic intensity, matrix-free implementations are well positioned to exploit this capability.

For quadrilaterals and hexahedra, classic nodal and modal bases with standard tensor-product structure admit efficient sum-factorisation, with many studies achieving remarkable performance on modern central processing units (CPUs) and graphics processing units (GPUs)~\cite{fischer2022nekrs,jansson2024neko,andrej2024mfem}. 
However, generating high-quality hexahedral meshes for complex geometries remains
challenging~\citep{shepherd2008HexConstraints,pietroni2022HexSurvey}.
Extending these techniques to three-dimensional unstructured meshes is therefore appealing
as it enables simulations in more complex industrial geometries.
This requires well-defined polynomial expansions on elements of general shapes,
including tetrahedra, prisms, and pyramids. Collapsed-coordinate constructions provide a generalised tensor-product
structure and corresponding orthogonal or modified modal bases~\citep{karniadakis_spectralhp_1999},
which underpin continuous and discontinuous Galerkin spectral/$hp$ methods on
unstructured meshes~\citep{warburton_galerkin_1999}; however, their efficient implementation has received less attention.
The open-source spectral element library Nektar++~\cite{cantwell2015,moxeyNektar2020} provides
these expansion bases for general element shapes. Sum-factorisation remains
applicable, although the operation counts and attainable performance
depend on the element shape and basis type~\citep{moxey2020}.
Recent work includes architecture-aware operator implementations
for mixed-element meshes~\cite{xing2026ArchitectureAware}, and unified
matrix-free interface-flux evaluation for general discontinuous Galerkin
spectral elements~\citep{xia2026}.

In real-world applications of discontinuous Galerkin methods, implicit solvers 
are often required, for which an efficient matrix-free operator is critical
for good performance, but not sufficient: an equally effective and
inexpensive preconditioner is also required, which remains a
challenge for DG spectral element methods.
A central issue is preconditioning for the Poisson equation, which is a recurring
task in the simulation of incompressible flow by splitting or projection
methods~\citep{karniadakis1991Splitting,guermond2006Projection}.
Unlike the positive mass-shifted Helmholtz systems arising in implicit viscous steps, 
the Poisson system lacks a mass term to make it more diagonally dominant.
Hybrid multigrid~\cite{fehn2020hybridMG} 
and low-order-refined (LOR) auxiliary-space preconditioning~\cite{pazner2020EfficientLOR}
have received substantial attention as promising solutions to this issue.

\citet{fehn2020hybridMG} developed a comprehensive
hybrid multigrid framework for symmetric interior-penalty Galerkin (SIPG) Poisson
problems on hexahedral nodal spectral elements,
combining discontinuous-to-continuous transfer ($c$), polynomial
($p$) and geometric ($h$) coarsening, together with algebraic multigrid (AMG) on the
coarsest grid. Chebyshev-accelerated point Jacobi is employed for smoothing, and
all components aside from the coarse-level solver are matrix-free.
Transferring to the continuous space at the finest level improves
robustness with respect to the penalty parameter in their experiments.
Hybrid coarsening also alleviated coarse-solver bottlenecks 
present in pure $p$- or $h$-multigrid. 
They also identified that a matrix-free smoother robust to anisotropic meshes
is an important direction for further improvement.
LOR preconditioning exploits spectral equivalence between
high-order operators and sparse low-order discretisations on refined
meshes. \citet{pazner2020EfficientLOR} combined finite-element/spectral-element (FEM--SEM)
equivalence on Gauss--Lobatto-refined meshes with overlapping additive
Schwarz methods, while retaining sum-factorised evaluation
of the high-order operator. Multigrid can serve as the auxiliary-space solver,
provided that it handles the increasing anisotropy of the refined meshes as $p$ increases.
For collapsed coordinates, \citet{pazner2026} used the Duffy transformation to
construct high-order spaces on triangles
and degree-independent spectral equivalence to low-order discretisations
on refined meshes.
These approaches avoid assembling the high-order matrix but still require an
assembled auxiliary-space system.

The cost of smoothing in multigrid has received particular attention.
Point Jacobi may be insufficient when the basis induces strong coupling within
elements. Block Jacobi, a non-overlapping additive Schwarz method,
resolves this coupling more effectively,
but storing and applying a dense elemental inverse becomes expensive at high order.
For modal DG discretisations, \citet{thiele2022} combined $p$-multigrid with AMG
and introduced partial smoothers that update only a subset of the unknowns on
each level. This reduces smoothing work while allowing a trade-off between
convergence and computational cost.

Several approaches reduce the cost of the local solves themselves.
\citet{bastian2019} replaced dense element-block factorisations by inexact,
matrix-free local Krylov solves; relatively loose local tolerances could preserve
outer convergence in their Poisson and diffusion experiments.
\citet{witteFastTensorProduct2021} used fast
diagonalisation for Cartesian cell and vertex-patch Schwarz smoothers,
with Cartesian surrogate operators extending the approach to non-Cartesian meshes.
Basis transformations provide another route to inexpensive local solves.
\citet{brubeck2022VertexStar} constructed a tensor-product basis that makes
separable vertex-patch operators sparse enough for efficient sparse Cholesky
factorisation. Nonseparable problems are treated using separable surrogate
operators. Their approach reduces the cost of factorising transformed patch
operators, rather than fitting an elemental inverse.

Algebraic tensor-product approximations also avoid the cost of full dense
elemental solves. \citet{pazner2018} constructed Kronecker-product approximations
of elemental Jacobian blocks and applied their inverses efficiently.
They also extend their method to symmetric interior-penalty discretisations~\citep{pazner2018IPTensor}.
\citet{diosady2019TensorProduct} approximated elemental Jacobian blocks by summations
of tensor products of small one-dimensional matrices, using an algebraic
procedure to minimise the approximation error. Formation, storage and
application scale linearly with solution order per degree of freedom in their
scalar formulation. These approaches approximate the elemental operator before
inversion and exploit standard tensor-product structure.

Other developments address implementation efficiency and element generality.
\citet{baccasImplementationTechniquesMultigrid2025} reduced data movement and
synchronisation in matrix-free multigrid, while
\citet{stillMatrixFreeEvaluationStrategies2025} demonstrated matrix-free
operators and hybrid multigrid on tetrahedral meshes at low-to-moderate
polynomial degrees.

Sparse approximate inverses offer a further way to reduce preconditioning costs.
They approximate the inverse, or a factor of it, by sparse matrices, so that
applying the preconditioner requires only sparse matrix--vector products.
Classical constructions either minimise a Frobenius-norm residual over an
adaptively selected sparsity pattern~\citep{grote1997}, or approximate the
inverse Cholesky factor, which preserves symmetry and positive
definiteness~\citep{kolotilina1993}.
For a fixed sparsity pattern, \citet{brarda2026SpectralSAI} instead select the
entries of a factorised approximate inverse by optimising spectral objectives that
cluster the eigenvalues of the preconditioned operator near unity, rather than
by minimising a Frobenius-norm residual.
Sparse approximate inverses have also been used as multigrid smoothers.
\citet{he2023OptimizedSAI} used local Fourier analysis to optimise such smoothers
for Laplacian systems, while \citet{jimenezrecio2026BlockFSAI} combined Chebyshev
smoothing with adaptive block-factorised approximate inverses for biharmonic and
triharmonic problems.
\citet{saye2025} compared block-sparse
approximate-inverse smoothers with block Jacobi and block Gauss--Seidel for local
discontinuous Galerkin discretisations of Poisson, elliptic interface and
multiphase Stokes problems. Their sparsity patterns range from single-element
blocks to multielement stencils, with entries obtained by minimising a
Frobenius-norm residual.

These studies motivate a balance between the accuracy and application cost of
local inverses. Extending standard tensor-product local-solver constructions to
collapsed-coordinate modal elements remains challenging because the basis
functions have coupled index ranges. We therefore seek an efficient approximate
elemental inverse with lower application cost that can compensate for additional
solver iterations.

Here we present a design that
constructs an approximate elemental inverse from inverse modal-to-nodal
transforms and an element-dependent sparse matrix. The transforms are
shared by elements of the same shape, basis, order and projection grid,
while the sparse matrix represents local operator and geometric information
in point space. Its entries are obtained by truncating a projected exact inverse
or by solving a weighted least-squares problem for the complete
transformed action. More importantly, standard, partial and generalised sum-factorisation
provide efficient inverse transforms for various shapes.
A related idea is used by \citet{ChanWarburton2017}, where they use the inverse 
of a reference mass matrix to approximate the exact inverse mass matrix.
We use the approximate elemental inverse as both a
standalone preconditioner and a multigrid Chebyshev smoother. Local spectral
analysis explains what the fitting objective controls and why symmetry
alone does not guarantee positivity. Benchmarks then compare convergence,
throughput and solve time with the elemental block inverse. 
The largest gains are observed on high-order hexahedra and prisms; tetrahedra and
pyramids show smaller, configuration-dependent benefits, compared with the exact element-block inverse.

The paper is organised as follows:
Section~\ref{sec:discretisation} introduces the SIPG formulation.
Section~\ref{sec:mg} describes the polynomial multigrid framework and smoother
choices, motivating an economical approximate elemental inverse.
Section~\ref{sec:approx-inverse} develops the modal-to-nodal transform and sparse
approximation, and then presents sum-factorised inverse
transforms and least-squares fitting of the sparse matrix entries.
Section~\ref{sec:results} presents the benchmark setup, the tuning of the
approximate elemental inverse as a standalone preconditioner and multigrid smoother, 
and the comparison with a continuous Galerkin solver on various workloads.

\section{Poisson problem and SIPG formulation}
\label{sec:discretisation}
We consider the Poisson problem on the domain $\Omega\subset\mathbb{R}^d$ with boundary $\partial\Omega=\Gamma=\Gamma_D\cup\Gamma_N$:
\begin{equation}
\begin{aligned}
    -\nabla^2 u &= f
    && \text{in } \Omega, \\
    u &= g
    && \text{on } \Gamma_D, \\
    \nabla u \cdot \mathbf{n} &= h
    && \text{on } \Gamma_N, 
\end{aligned}
 \label{eq:poisson}
\end{equation}
Here $d$ is the spatial dimension, $f$ is the source term, $g$ and $h$
are the prescribed Dirichlet and Neumann data, and $\boldsymbol{n}$ is the
outward unit normal.
The computational domain $\Omega_h$ is partitioned into elements $\Omega^e$, each mapped from a
reference element. The numerical solution $u_h$ in the discontinuous spectral
element space $\mathcal V_h$ of order $p$ is defined by a shape-dependent modal expansion on the reference
element $\Omega_{\text{st}}$
\begin{equation} \label{eq:expansion}
  u_h(\boldsymbol{\xi})=\sum_{i}{\phi_i(\boldsymbol{\xi})\hat{u}_i}~,
\end{equation}
where $\hat{u}_i$ is the coefficient associated with the $i$th basis function. 
By collecting all coefficients of an element into a vector, 
we can rewrite the above formulation in matrix form 
\begin{equation} \label{eq:bwdtrans}
    \mathbf{u}_h=\mathbf{B}\widehat{\mathbf{u}}, 
\end{equation}
where $\mathbf{u}_h$ contains values at evaluation points $\boldsymbol{\xi}_j$
and $\mathbf{B}$ is the basis matrix with entries $B_{ji}=\phi_i(\boldsymbol{\xi}_j)$.
Detailed formulations of the modal expansion bases $\phi_i$ for various element shapes
can be found in \citet{karniadakis_spectralhp_1999}. These modal expansion bases
have either standard tensor-product structures (on quadrilateral and hexahedral elements)
or generalised tensor-product structures (on triangular, tetrahedral, prismatic and pyramidal elements),
which allow sum-factorisation to reduce the computational cost and memory traffic \cite{moxey2020}.
The local coordinate $\boldsymbol{\xi}\in\Omega_{\text{st}}$
is linked with the physical coordinate $\boldsymbol{x}\in\Omega^e$ by a geometric mapping
\begin{equation}
\boldsymbol{x}=\boldsymbol{\chi}^e(\boldsymbol{\xi}),
\label{eq:geometry-map}
\end{equation}
from which the inverse Jacobian matrix $\mathbf{J}=\frac{\partial \boldsymbol{\xi}}{\partial \boldsymbol{x}}$
and its determinant $J$ are derived and used for finite element integration and differentiation.

For discontinuous discretisations, the data on either side of an element interface may differ.
We use superscripts $^-$ and $^+$ to denote quantities or objects from the current and neighbouring
elements. The average and jump are defined by
\begin{equation}
 \avg{u_h}=\frac{u_h^-+u_h^+}{2}~,\qquad
 \jump{u_h}=u_h^-\boldsymbol{n}^-+u_h^+\boldsymbol{n}^+~,
\end{equation}
where $\boldsymbol{n}^{\pm}$ are the outward unit normals of the adjacent
elements $\Omega^{\pm}$. Vector averages are defined componentwise. These quantities are used
to describe the numerical flux across element boundaries.
$(\cdot,\cdot)_{\Omega^e}$
and $(\cdot,\cdot)_{\partial\Omega^e}$ denote volume and surface inner products, respectively.
We seek $u_h\in\mathcal V_h$ such that, for every local test function $v_h$,
\begin{equation}
 \begin{aligned}
 a_h^e(v_h,u_h)
 ={}&(\nabla v_h,\nabla u_h)_{\Omega^e}
 -(\nabla v_h,\tfrac12\jump{u_h})_{\partial\Omega^e}\\
 &-(v_h,\avg{\nabla u_h}\cdot\boldsymbol{n})_{\partial\Omega^e}
 +(v_h,\tau\jump{u_h}\cdot\boldsymbol{n})_{\partial\Omega^e}
 =(v_h,f)_{\Omega^e},
 \end{aligned}
 \label{eq:sipg}
\end{equation}
where $\boldsymbol{n}$ is short for $\boldsymbol{n}^-$. This is the Poisson specialisation
of the symmetric interior-penalty Galerkin (SIPG) formulation on a single element \cite{arnold2002, xia2026}.
The three surface terms are the symmetric consistency correction, 
average-gradient flux and jump penalty.

Boundary conditions are imposed weakly through the exterior values. 
For example, on the Dirichlet boundary $\Gamma_D$,
\begin{equation}
 u_h^+=2g-u_h^-,\qquad
 \nabla u_h^+\cdot\boldsymbol{n}=\nabla u_h^-\cdot\boldsymbol{n}
~.
 \label{eq:dirichlet}
\end{equation}
On the Neumann boundary $\Gamma_N$,
\begin{equation}
 u_h^+=u_h^-,\qquad
 \nabla u_h^+\cdot\boldsymbol{n}=2h - \nabla u_h^-\cdot\boldsymbol{n} ~.
 \label{eq:neumann}
\end{equation}
The penalty $\tau$ depends on order and element geometry through
a local trace estimate \citep{shahbazi2005,hillewaert2013development}
\begin{equation} \label{eq:penalty-parameter}
    \tau= \frac{\sigma_0}{V(\Omega^e)}\left(\frac{1}{2}\sum_{f\in \partial\Omega^e \backslash\Gamma}C_{e,f}A(f)
        + \sum_{f\in \partial\Omega^e \cap\Gamma}C_{e,f}A(f)\right)
\end{equation}
where $V(\Omega^e)$ is the element volume, $A(f)$ is the area of face $f$,
and $C_{e,f}$ is a dimensionless trace-estimate coefficient depending on the
polynomial degree, element shape and face type.
\citet{warburton2003TraceInverse} derived these coefficients explicitly for
triangles and tetrahedra; we use the shape- and face-dependent values given by
\citet{hillewaert2013development}.
We additionally scale the penalty by a dimensionless factor $\sigma_0$ 
to tune the numerical stability across different cases. 
We always choose the biggest $\tau=\text{max}(\tau^{-},\tau^{+})$ 
across the element interfaces to ensure a symmetric and stable discretisation.

After transferring prescribed-boundary contributions to the right-hand
side, the discrete problem $\mathbf{A}\widehat{\mathbf{u}}=\mathbf{b}$ is formed. 
Here $\widehat{\mathbf{u}}$ collects the unknown modal coefficients, 
and $\mathbf{b}$ contains the forcing and boundary
contributions. Under the assumptions above, the global SIPG matrix $\mathbf{A}$
is symmetric positive definite (SPD), so the discrete system can be solved
by the conjugate gradient method.
The convergence depends on the conditioning of the original system, as well as the
effectiveness of the preconditioner.

The SIPG formulation and Poisson solver are implemented in the open-source spectral element framework
Nektar++. Matrix-free and sum-factorisation techniques are employed in the evaluation of the SIPG formulation 
to reduce the memory footprint and improve cache locality and performance.
Single-instruction, multiple-data (SIMD) instructions are effectively utilised to boost the floating-point
performance on CPUs, thanks to the cache-friendly data layout.
Detailed performance characteristics and the generalisation of the matrix-free DG framework 
to different shapes of spectral elements can be found in~\citep{xia2026}.

\section{Polynomial multigrid preconditioning}
\label{sec:mg}
For unstructured high-order finite element discretisations,
polynomial multigrid ($p$-multigrid, or $p$-MG) is a convenient choice because
it avoids mesh agglomeration.
\citet{fehn2020hybridMG} provide a reference framework for $p$-MG.
The standard preconditioned conjugate gradient method requires a fixed, linear,
SPD preconditioner. Violating these assumptions removes its usual convergence
guarantees.
The structure of the $p$-MG used in this work is described as follows.

\subsection{Coarse-level solver}
In $p$-multigrid, we always choose the lowest order $p=1$ as the coarsest level. 
Since we exclude geometric multigrid in this work, the natural choice is to use algebraic multigrid (AMG)
for this low-order problem. Applying a single iteration of AMG is not strong enough to reduce the residual and obtain
an effective approximation. \citet{fehn2020hybridMG} use an inner Krylov solver on the coarse level
and use AMG as its preconditioner. The inner Krylov
solve can make the preconditioner nonstationary and may require a flexible
outer solver, although they observed no convergence problems.
We use Chebyshev iteration with an AMG preconditioner as the coarse level solver.
With a fixed iteration count, spectral bounds and a fixed linear AMG
preconditioner, this defines a fixed linear and SPD operation.
Whatever the choice, the system matrix on the coarse level needs to be assembled, and the cost
is acceptable as the low-order problem is much smaller than the original high-order one.

\subsection{Restriction and prolongation} 
The prolongation operator $\mathcal{P}$ maps coarse-level coefficients at
level $l$ to the finer level $l-1$. We use its transpose
$\mathcal{R}=\mathcal{P}^\T$ to restrict residuals.
An $L_2$-projection construction can be written as follows:
\begin{equation}
\begin{aligned}
    \widehat{\mathbf{u}}^{l-1} =\mathcal{P}\widehat{\mathbf{u}}^{l} &=  (\mathbf{M}^{l-1})^{-1}~(\mathbf{B}^{l-1})^\T~ \mathbf{W}~\mathbf{B}^{l}~\widehat{\mathbf{u}}^{l} \\
    \mathbf{r}^{l} =\mathcal{R}\mathbf{r}^{l-1} &=  (\mathbf{B}^{l})^\T~\mathbf{W}~\mathbf{B}^{l-1}~(\mathbf{M}^{l-1})^{-1}~\mathbf{r}^{l-1}~,
\end{aligned}
\end{equation}
where $\mathbf{r}^{l}$ is the residual on level $l$. The matrices
$\mathbf{B}^{l}$ and $\mathbf{B}^{l-1}$ evaluate the respective bases on a
common quadrature grid. The mass matrix is $\mathbf{M}=\mathbf{B}^\T\mathbf{W}\mathbf{B}$,
where $\mathbf{W}$ contains quadrature weights and Jacobians $\omega J$.
A simpler construction uses basis interpolation:
\begin{equation}
\begin{aligned}
    \widehat{\mathbf{u}}^{l-1} =\mathcal{P}\widehat{\mathbf{u}}^{l} &=  (\mathbf{B}^{l-1})^{+}~\mathbf{B}^{l}~\widehat{\mathbf{u}}^{l} \\
    \mathbf{r}^{l} =\mathcal{R}\mathbf{r}^{l-1} &=  (\mathbf{B}^{l})^\T~[(\mathbf{B}^{l-1})^{+}]^\T ~\mathbf{r}^{l-1}~,
\end{aligned}
\end{equation}
where the Moore--Penrose pseudoinverse $\mathbf{B}^+$ can be used when $\mathbf{B}$ is rectangular, 
particularly in collapsed-coordinate expansions. 
For nested spaces and full-column-rank basis matrices, this gives the same
prolongation as $L_2$ projection because coarse function spaces are included
in the fine space. This formulation avoids the explicit storage of the inverse mass matrix for each element 
and makes $\mathcal{P}$ and $\mathcal{R}$ shared matrices for all elements of the same type, 
which can be implemented by general matrix--matrix multiply (GEMM).

\subsection{Smoother choices}
Smoothing is critical in the multigrid method. The solution projected from the coarse level is still an approximation and requires smoothing to reduce the remaining error. 
From another perspective, the low-rank coarse solver cannot be used as a preconditioner for a higher-rank solver directly due to rank deficiency.
Existing convergence analyses require the number of smoothing steps to grow with
the polynomial degree.
For $hp$-DG discretisations of elliptic problems, \citet{antonietti2015} proved
level- and discretisation-independent convergence of W-cycle $hp$-multigrid,
provided that the number of smoothing steps grows as $O(p^2)$ or $O(p^3)$.
For SIPG with hierarchical Legendre bases, \citet{lei2025} reduced this
requirement to $O(p)$ for $p$-multigrid.
The cost of each smoothing step therefore becomes increasingly important at high order.
Polynomial smoothers offer computational advantages for parallel multigrid on
SPD systems~\citep{adams2003PolynomialSmoothing}.
Following \citet{fehn2020hybridMG}, we employ Chebyshev iterations as both
pre- and post-smoothers to achieve a higher convergence rate on symmetric systems. The effectiveness and cost strongly depend on the local
preconditioner used within Chebyshev smoothing. AMG is not a good option here, since the smoother
will be used at the finest level and assembling the system matrix for high-order problems is extremely expensive. 
Two common choices are point-Jacobi and block-Jacobi preconditioners:
\begin{equation}
 \begin{aligned}
 \text{Point Jacobi:}\quad &{z}_i=(A_{ii})^{-1}{r}_i~,\\
 \text{Block Jacobi:}\quad &\mathbf{z}_e=\mathbf{H}_e\mathbf{r}_e~,\qquad \mathbf{H}_e=\mathbf{A}_e^{-1},
 \end{aligned}
 \label{eq:local-preconditioners}
\end{equation}
where the subscript $e$ denotes the part of the vector or matrix that corresponds to element $e$.
$\mathbf{A}_e$ is the elemental diagonal block of $\mathbf{A}$.
Both point Jacobi and block Jacobi act independently on each element. 
Point Jacobi rescales each coefficient separately; 
block Jacobi retains all coupling within $\mathbf{A}_e$
and is equivalent to solving $\mathbf{A}_e\mathbf{z}_e=\mathbf{r}_e$.
Point Jacobi is inexpensive, but cannot provide sufficient smoothing for
strongly coupled modal bases. Block Jacobi treats these couplings more
effectively, but storing and applying the dense inverse $\mathbf{H}_e$ is costly at
high order. We therefore seek an approximation of the exact inverse, $\mathbf{P}_e \approx \mathbf{H}_e$, 
along with a cheap factorisation of $\mathbf{P}_e$, to replace the dense
matrix--vector multiplication in block Jacobi. We call $\mathbf{P}_e$
the \emph{approximate elemental inverse} or \emph{approximate inverse} 
as a short form when the context is clear.

\section{Design of the approximate elemental inverse}
\label{sec:approx-inverse}
\subsection{Modal-to-nodal transform and sparse inverse approximation}
\label{sec:construction}

From the computational perspective, an ideal structure of the preconditioner is
\begin{equation}
 \mathbf{P}_e=\mathbf{L}\mathbf{V}_e\mathbf{L}^\T\approx \mathbf{H}_e,
 \label{eq:approx}
\end{equation}
where $\mathbf{L}$ depends only on the reference basis and is shared by 
elements of the same shape, basis, polynomial degree and projection grid.
The symmetric sparse matrix $\mathbf{V}_e$ contains the element-dependent
geometry and local operator information. The repeated general matrix--vector multiply (GEMV) 
of $\mathbf{L}$ on different element vectors can be batched into GEMM, 
increasing reuse of the shared basis data and offering higher arithmetic intensity. 
Sparsity in $\mathbf{V}_e$ reduces element-dependent storage, memory traffic and multiplication cost.
Efficient factorisations of $\mathbf{L}$, described in Section~\ref{sec:inverse-sum-factorisation}, further
reduce the transform work.

We construct $\mathbf{L}$ from an inverse modal-to-nodal
transform.
This is motivated by the favourable conditioning and
iterative convergence of tensor-product Gauss--Lobatto--Legendre (GLL)
nodal bases on hexahedral elements. Compared with modal bases,
the corresponding elemental blocks exhibit a clearer sparsity pattern,
with dominant entries more concentrated near the diagonal under a
tensor-product ordering. Moreover, nodal degrees
of freedom are associated with spatial locations, making proximity
between points a natural guide for selecting a sparse pattern of $\mathbf{V}_e$.

\subsubsection{Projection space}
\label{sec:projection}
Here, a \emph{projection space} denotes the representation of an elemental field
by its values on a selected evaluation grid.
For example, tensor-product GLL grids provide a nodal representation on hexahedral elements.
For general element shapes, we consider the following projection spaces:
\begin{itemize}
    \item \textbf{Tensorial GLL grid.} This is the standard nodal grid for hexahedral
    or quadrilateral elements, using GLL points.
    We use $p+1$ GLL points in each coordinate direction, giving the same total
    number of points as modal coefficients on hexahedral and quadrilateral elements.
    For collapsed-coordinate elements,
    the total number of points exceeds the number of coefficients,
    and multiple points can coincide at collapsed vertices or edges. 
    We denote the corresponding projection as \texttt{ModalToGLL} or \texttt{ToGLL} for short in this work.
    \item \textbf{Original quadrature grid.} This uses tensor-product
    quadrature points for inner products within the element.
    For degree $p$, we typically use $p+1$ Gaussian points 
    or $p+2$ Gauss--Lobatto points on each coordinate. 
    $p+1$ Gauss--Radau points are used on collapsed coordinates to avoid the top-vertex singularity issue \cite{karniadakis_spectralhp_1999}. 
    The total point count again exceeds the coefficient count in collapsed-coordinate elements.
    The field representation based on the quadrature grid is termed \textit{physical space} in Nektar++~\cite{xia2026}, 
    so we denote the corresponding projection as \texttt{ModalToPhys} or \texttt{ToPhys} in this work.
    \item \textbf{Non-tensorial nodal point set.} This uses rotationally
    symmetric point sets, such as electrostatic points on
    triangles~\citep{hesthaven1998Electrostatics} or tetrahedra~\citep{hesthaven2000TetNodes}.
    For prisms, we can use GLL points in the extrusion
    direction, paired with triangular electrostatic points.
    The total number of points matches the total number of coefficients, which can potentially
    save computation and storage. We denote the corresponding projection as \texttt{ModalToNonTensorNodal} or \texttt{ToNonTensorNodal} for short.
\end{itemize}
Let $N_P$ be the total number of modal coefficients and $N_Q$ the total
number of points on the selected evaluation grid. The basis matrix
$\mathbf{B}\in\R^{N_Q\times N_P}$ is assumed to have full column rank
for the selected point set; $N_Q\geq N_P$ is necessary but not sufficient.
The task is to find $\mathbf{L}\in\R^{N_P\times N_Q}$ such that
\begin{equation}
 \mathbf{L}\mathbf{B}=\mathbf{I}_{N_P}.
 \label{eq:leftinverse}
\end{equation}
Here $\mathbf{I}_{N_P}$ is the $N_P\times N_P$ identity matrix.
If $N_Q=N_P$, $\mathbf{L}$ is the exact inverse. If $N_Q>N_P$, $\mathbf{L}$ can be the Moore--Penrose pseudoinverse
$\mathbf{B}^+$. Then $\mathbf{L}\mathbf{u}$ minimises
$\|\mathbf{B}\widehat{\mathbf{u}}-\mathbf{u}\|_2$, and
$\mathbf{L}^\T\widehat{\mathbf{u}}$ is the minimum-norm solution of
$\mathbf{B}^\T\mathbf{x}=\widehat{\mathbf{u}}$. Other left inverses satisfy
Eq.~\eqref{eq:leftinverse} without these optimality properties.

\subsubsection{Sparse patterns}
\label{sec:sparse-pattern}
Using a sparse matrix $\mathbf{V}_e$ is critical to achieving lower computational cost
and higher efficiency. However, a sparser pattern reduces storage and computation
but can also discard couplings needed to approximate the inverse of the elemental block accurately. The
choices of sparse pattern need to be tested and compared.
In this work, we tested the following sparse patterns:
\begin{itemize}
 \item \texttt{Band<N>} \quad A banded matrix with $N$ non-zero diagonals,
 where $N\geq1$ is an odd number. \texttt{Band1} can also be called 
 \texttt{Diagonal}.
 \item \texttt{NearestPnt} \quad  The nearest-point stencil on a structured
 grid, with up to 5 entries per row in two dimensions and 7 in
 three dimensions, including the diagonal.
 \item \texttt{FirstLayer} \quad  The first-layer stencil on a structured grid,
 with up to 9 entries per row in two dimensions and 27 in three
 dimensions.
\end{itemize}
The nearest-point and first-layer stencils are designed for the tensorial point sets,
with the assumption that the coupling between nearest points should not be ignored.
The diagonal matrix requires the least storage, while the first-layer stencil requires the most storage.
In practice, $\mathbf{V}_e$ needs to be sparse enough to compensate for the additional 
computing cost of transform operations, so we do not use a high-bandwidth sparse matrix.

Let $\mathcal S$ be a symmetric sparse pattern, and let operator $\mathcal M_{\mathcal S}$ 
retain supported entries and set the others to zero. 
The simplest construction of $\mathbf{V}_e$ truncates the projected inverse:
\begin{equation}
 \mathbf{V}_e^{\rm proj}=\mathcal M_{\mathcal S}(\mathbf{B}\mathbf{H}_e\mathbf{B}^\T).
 \label{eq:projected}
\end{equation}

Without truncation, every exact left inverse recovers $\mathbf{H}_e$:
\begin{equation}
 \mathbf{L}(\mathbf{B}\mathbf{H}_e\mathbf{B}^\T)\mathbf{L}^\T=\mathbf{H}_e.
 \label{eq:fullreconstruction}
\end{equation}
With truncation, the defect of the projected values is
\begin{equation}
 \begin{aligned}
 \mathbf{\Delta}_e&=\mathcal M_{\mathcal S}(\mathbf{B}\mathbf{H}_e\mathbf{B}^\T)-\mathbf{B}\mathbf{H}_e\mathbf{B}^\T.
 \end{aligned}
 \label{eq:defect}
\end{equation}
Thus different left inverses can affect the quality of approximation, even though we use the same sparse pattern.
We next describe efficient implementations of $\mathbf{L}$, then optionally
use least-squares fitting to optimise the values of the sparse matrix.

\subsection{Sum-factorisation for inverse transforms}
\label{sec:inverse-sum-factorisation}
Storing a dense representation of the left inverse $\mathbf{L}$ and performing GEMM during preconditioning
would cause substantial transform work that is far more expensive than applying $\mathbf{V}_e$.
The tensorial structure of the nodal space and modal expansion basis gives us the
opportunity to exploit sum-factorisation when applying the inverse of the basis matrix,
and significantly reduce the computational cost.

\subsubsection{Tensorial and partial tensorial sum-factorisation}
For hexahedra, the basis functions are products of three one-dimensional
basis functions. If the evaluation points also form a tensor-product grid, then the basis matrix
can be written as
\begin{equation}
 \mathbf{B}=\mathbf{B}_3\otimes \mathbf{B}_2\otimes \mathbf{B}_1,
 \label{eq:tensor}
\end{equation}
where $\otimes$ denotes the Kronecker product and $\mathbf{B}_d$ is the one-dimensional basis matrix that evaluates
the basis values on the given points in direction $d$.
The backward transform (\texttt{BwdTrans}), defined in Eq.~\eqref{eq:bwdtrans},
is then evaluated by three successive one-dimensional contractions, 
giving standard sum-factorisation \cite{karniadakis_spectralhp_1999}
\begin{equation}
 (\mathbf{B}\widehat{\mathbf{u}})_{rst}
 =\sum_{i=0}^{p}\phi_i(\xi_{1r})
 \left\{\sum_{j=0}^{p}\phi_j(\xi_{2s})
 \left[\sum_{k=0}^{p}\widehat{u}_{ijk}\phi_k(\xi_{3t})\right]\right\},
 \label{eq:hex-sum-factorisation}
\end{equation}
where $r$, $s$ and $t$ index the evaluation points in the three directions.
When the basis matrices are square and invertible, 
the inverse has the same tensor-product
structure:
\begin{equation}
 \mathbf{L}=\mathbf{B}^{-1}=\mathbf{B}_3^{-1}\otimes \mathbf{B}_2^{-1}\otimes \mathbf{B}_1^{-1}.
 \label{eq:tensorinverse}
\end{equation}
We denote this construction by \texttt{SumFacInverse}. For rectangular
basis matrices, the Moore--Penrose pseudoinverse can be factorised similarly;
we denote this construction by \texttt{SumFacPseudoInv}:
\begin{equation}
 \mathbf{L}=\mathbf{B}^+=\mathbf{B}_3^+\otimes \mathbf{B}_2^+\otimes \mathbf{B}_1^+ .
 \label{eq:tensorpseudoinverse}
\end{equation}

For prisms, the basis can be considered as a tensor product of a two-dimensional triangular
basis and a one-dimensional basis for the extrusion direction, which results in
a partial sum-factorisation like
\begin{equation}
 \mathbf{B}=\mathbf{B}_{\rm 1D}\otimes \mathbf{B}_{\triangle},\qquad
 \mathbf{L}=\mathbf{B}^+=\mathbf{B}_{\rm 1D}^+\otimes \mathbf{B}_{\triangle}^+,
 \label{eq:prismtensor}
\end{equation}
where $\mathbf{B}_{\triangle}$ denotes the full basis matrix of the triangular expansion.
The pseudoinverse can still be recovered in this way.

\subsubsection{Generalised sum-factorisation in collapsed coordinates}
\label{sec:collapsed}
For tetrahedra and pyramids, neither tensorial nor partial tensorial factorisation exists.
However, the \texttt{BwdTrans} matrix in collapsed-coordinate elements can still be 
factorised by the so-called generalised sum-factorisation \cite{karniadakis_spectralhp_1999,moxey2020}. 
Taking the modified modal basis on a tetrahedron as an example, 
$\mathbf{B}$ denotes the full \texttt{BwdTrans} matrix
and $\mathbf{B}'$ denotes the matrix without the additional correction. 
The correction on vertex and edge terms forms a matrix denoted as
$\Delta \mathbf{B}$, so that
\begin{equation}
 \mathbf{B}=\mathbf{B}'+\Delta \mathbf{B}.
 \label{eq:basiscorrection}
\end{equation}
All three matrices have $N_Q$ rows and $N_P$ columns. 
For $\Delta \mathbf{B}$, only the columns that correspond to the modes that need to be corrected
are non-zero.

Let $\phi_i^a$, $\phi_{ij}^b$ and $\phi_{ijk}^c$ denote the Jacobi-based
functions in the three collapsed coordinates $\boldsymbol{\eta}$. 
For a degree-$p$ expansion, the expansion value at the point $(r,s,t)$ in a tensorial nodal grid
can be evaluated by
\begin{equation}
 u_q=(\mathbf{B}'\widehat{\mathbf{u}})_q
 =\sum_{i=0}^{p}\phi_i^a(\eta_{1r})
 \left\{\sum_{j=0}^{p-i}\phi_{ij}^b(\eta_{2s})
 \left[\sum_{k=0}^{p-i-j}\widehat u_{ijk}
 \phi_{ijk}^c(\eta_{3t})\right]\right\}.
 \label{eq:tetsumfac}
\end{equation}
Here $q$ is the flattened index corresponding to $(r,s,t)$.
We first contract over $k$, then over $j$, and finally over $i$.
The dependence of $\phi_{ij}^b$ on $i$ and of $\phi_{ijk}^c$ on $(i,j)$
indicates this is not the Kronecker product of three independent factors.

Although reversing these stages does not generally recover the pseudoinverse,
we can still construct a left inverse $\mathbf{L}'$ of $\mathbf{B}'$ 
by reversing the contraction stages given in Eq.~\eqref{eq:tetsumfac}. 
Each stage uses the pseudoinverse of its conditional one-dimensional
basis matrix that depends on the previous indices $i$ or $(i,j)$. 
If all these matrices have full column rank, the resulting map satisfies
\begin{equation}
 \mathbf{L}'\mathbf{B}'=\mathbf{I}_{N_P}.
 \label{eq:productleftinverse}
\end{equation}
Since $\mathbf{L}'\mathbf{B}=\mathbf{I}_{N_P}+\mathbf{L}'\Delta \mathbf{B}$, 
the additional correction terms are then included via
\begin{equation}
 \mathbf{L}=\left(\mathbf{I}_{N_P}+\mathbf{L}'\Delta \mathbf{B}\right)^{-1}\mathbf{L}',
 \label{eq:correction}
\end{equation}
provided the correction matrix is nonsingular. This gives
$\mathbf{L}\mathbf{B}=\mathbf{I}_{N_P}$.
The limited non-zero columns in $\Delta \mathbf{B}$ allow this correction matrix
to be constructed as a small sparse system during setup.
During the inverse transform, we first evaluate $\mathbf{L}'\mathbf{u}$, then apply the correction.
The resulting \texttt{SumFacLeftInv} is an exact left inverse of $\mathbf{B}$,
but generally differs from $\mathbf{B}^+$.

The same principle extends to pyramids and prisms, 
which implies that for prisms, we can do both \texttt{SumFacPseudoInv} and \texttt{SumFacLeftInv}.
Apart from sum-factorisation, the full dense pseudoinverse, denoted as \texttt{FullPseudoInv},
is also benchmarked and compared with the sum-factorised left inverse.

\subsection{Implementation and performance optimisation}
\label{sec:cost}
Here matrix-free refers to the repeated preconditioner application.
The current setup still constructs dense elemental inverses to obtain
the projected or fitted sparse entries.
Let $P=p+1$ denoting the number of degrees of freedom or unknown coefficients in each direction.
Both $N_P$ and $N_Q$ then scale as $O(P^3)$ in three dimensions, with shape-dependent constants. 
For the inverse transform $\mathbf{L}$, a dense representation requires $N_PN_Q$ entries and $O(N_PN_Q)$
floating-point operations per application, both scaling as $O(P^6)$.
Sum-factorisation of $\mathbf{L}$ in hexahedra requires $O(P^2)$ storage and $O(P^4)$ floating-point operations.
The generalised sum-factorisation in tetrahedra and pyramids requires $O(P^3)$ storage and $O(P^4)$ floating-point operations.
The partial sum-factorisation in prisms requires $O(P^4)$ storage and $O(P^5)$ floating-point operations 
because its triangular basis matrix is a dense block. 
A detailed formulation of the theoretical floating-point operation counts of the sum-factorisation in various elements can be 
found in \citet{xing2026ArchitectureAware}. In short, sum-factorisation significantly reduces the floating-point operation counts
over the full dense matrix multiplication, particularly in high-order cases. However, the actual performance gains still 
depend on the specific implementation and coding practices.
The sum-factorisation of $\mathbf{L}$ is implemented following the design of the \texttt{BwdTrans} operator in \citet{moxey2020}. We interleave the vector data from 4 or 8 elements, pack the $i$-th entries of the 4 or 8 elements
(depending on SIMD width) into a SIMD register, and then perform SIMD arithmetic operations upon them. 
Compile-time constant loop counts allow the compiler to do more aggressive loop unrolling or fusion to optimise the pipelines. 
This design has been shown to exploit floating-point performance on CPUs effectively.
The full dense matrix multiplication of $\mathbf{L}$ is implemented by packing multiple element vectors into a matrix and forming a GEMM operation. 
However, applying \texttt{LIBXSMM} \cite{heineckeLIBXSMM2016} over small batches
does not perform better than applying OpenBLAS once over all elements, 
since the dense inverse transform is large and matrix tiling is more important than effective SIMD dispatch.

The sparse matrix $\mathbf{V}_e$ stores and accesses $O(P^3)$ entries per element. 
Element vectors and temporary workspace also require $O(P^3)$ storage. 
$\mathbf{V}_e$ is not shared by multiple elements; 
thus, it accounts for most of the element-dependent storage and memory traffic, and the arithmetic intensity can be low.
Since the application of $\mathbf{V}_e$ is between two inverse transforms, 
it is straightforward to use interleaved data and Single-instruction, multiple-data (SIMD) instructions 
to perform the sparse matrix--vector multiplication.
Pattern-specific kernels, rather than a general sparse matrix multiplication kernel, can be used
to avoid using additional sparse indices and improve the performance.

For the exact dense inverse of an elemental diagonal block, storage is
$N_P^2$ entries per element and application costs $O(N_P^2)$ floating-point
operations; both scale as $O(P^6)$. We cannot apply GEMM here, thus the potential performance speedup is limited.
The floating-point cost of applying two full dense inverse transforms 
is definitely much higher than applying an exact inverse of the elemental block; 
thus, sum-factorisation of the inverse transform is the key 
to making the approximate inverse more cost-effective than the exact inverse.

\subsection{Least-squares fitting of sparse matrix entries}
\label{sec:least-squares-fit}
\subsubsection{Energy-scaled least-squares problem}
With the transform $\mathbf{L}$ and the sparse pattern $\mathcal{S}$ fixed,
a natural idea is to determine the entries of $\mathbf{V}_e$ by
fitting the resulting approximate inverse $\mathbf{P}_e$ to $\mathbf{H}_e$.
Truncation in Eq.~\eqref{eq:projected} is straightforward but its final 
approximation quality strongly depends on the matrix structure after the transform
and the chosen sparse pattern.

A straightforward objective is to minimise the error between $\mathbf{P}_e$ and $\mathbf{H}_e$.
\begin{equation}
\|\mathbf{H}_e-\mathbf{P}_e\|_F~.
\end{equation}
Instead, we multiply the error matrix on both sides by $\mathbf{A}_e^{1/2}$, the symmetric positive-definite
square root of $\mathbf{A}_e$. Since $\mathbf{H}_e$ is the exact inverse of $\mathbf{A}_e$,
this scaling maps the exact inverse to the identity, and then we get the weighted Frobenius norm
\begin{equation}
 \|\mathbf{H}_e-\mathbf{P}_e\|_{F(\mathbf{A}_e)}
 :=\|\mathbf{A}_e^{1/2}(\mathbf{H}_e-\mathbf{P}_e)\mathbf{A}_e^{1/2}\|_F
 =\|\mathbf{I}_{N_P}-\mathbf{A}_e^{1/2}\mathbf{P}_e\mathbf{A}_e^{1/2}\|_F.
 \label{eq:weighted-frobenius-norm}
\end{equation}
To evaluate this objective without forming $\mathbf{A}_e^{1/2}$,
we define
\begin{equation}
 \mathbf{T}_e=\mathbf{C}_e^{-1}\mathbf{L},
 \label{eq:fit-transform}
\end{equation}
where $\mathbf{C}_e$ is the Cholesky factor of
$\mathbf{H}_e=\mathbf{C}_e\mathbf{C}_e^\T$. Then
\begin{equation*}
 \begin{aligned}
  \mathbf{A}_e^{1/2}\mathbf{P}_e\mathbf{A}_e^{1/2}
  &=\mathbf{A}_e^{1/2}\mathbf{L}\mathbf{V}_e
    \mathbf{L}^\T\mathbf{A}_e^{1/2}\\
  &=(\mathbf{A}_e^{1/2}\mathbf{C}_e)
    \mathbf{T}_e\mathbf{V}_e\mathbf{T}_e^\T
    (\mathbf{A}_e^{1/2}\mathbf{C}_e)^\T.
 \end{aligned}
\end{equation*}
Since $\mathbf{A}_e^{1/2}\mathbf{C}_e$ is orthogonal, this similarity
transformation preserves the Frobenius distance from the identity.
The weighted least-squares problem is therefore
\begin{equation}
\label{eq:fitting-problem}
\begin{aligned}
\min_{\mathbf{V}_e}\quad&
\frac{1}{2}
\left\|
\mathbf{I}_{N_P}
-\mathbf{T}_e\mathbf{V}_e\mathbf{T}_e^\top
\right\|_F^2,
\\
\text{subject to}\quad&
\mathbf{V}_e=\mathbf{V}_e^\top,
\qquad
\operatorname{supp}(\mathbf{V}_e)\subseteq\mathcal{S}.
\end{aligned}
\end{equation}
Here $\operatorname{supp}$ denotes the set of non-zero matrix positions.
This is a linear least-squares problem with the upper triangular part of $\mathbf{V}_e$ as independent unknowns. 
This error norm is related to eigenvalues of the preconditioned matrix, that is,
\begin{equation}
 \|\mathbf{I}_{N_P}-\mathbf{T}_e\mathbf{V}_e\mathbf{T}_e^\T\|_F^2
 =\|\mathbf{I}_{N_P}-\mathbf{A}_e^{1/2}\mathbf{P}_e\mathbf{A}_e^{1/2}\|_F^2
 =\sum_{j=1}^{N_P}(1-\lambda_j)^2,
 \label{eq:fitspectrum}
\end{equation}
where $\lambda_j$ are the eigenvalues of $\mathbf{A}_e^{1/2}\mathbf{P}_e\mathbf{A}_e^{1/2}$. 
In other words, the objective actually penalises the
deviation of the eigenvalues from unity. It does not
directly minimise their condition number or guarantee faster global
convergence. Its advantage over the unweighted norm is its direct
interpretation in terms of the preconditioned spectrum.

Setting the derivatives with respect to the independent entries of $\mathbf{V}_e$
to zero yields the normal equations
\begin{equation}
 \mathcal{M}_{\mathcal{S}}
 \left(\mathbf{T}_e^\top\mathbf{T}_e\mathbf{V}_e\mathbf{T}_e^\top\mathbf{T}_e\right)
 = \mathcal{M}_{\mathcal{S}}\left(\mathbf{T}_e^\top\mathbf{T}_e\right),
 \label{eq:fitnormal}
\end{equation}
where $\mathcal{M}_{\mathcal{S}}$ retains entries in the prescribed
symmetric sparsity pattern $\mathcal{S}$ and sets all other entries
to zero. These equations are solved to obtain independent entries of $\mathbf{V}_e$
using a diagonally preconditioned conjugate-gradient
least-squares method (CGLS). The left-hand side of \eqref{eq:fitnormal} is evaluated by
successive products with $\mathbf{T}_e$ and $\mathbf{T}_e^\top$,
followed by restriction to $\mathcal{S}$, without explicitly
assembling the least-squares design matrix or its normal matrix.
We use at most 30 iterations to obtain an approximate solution,
starting from the sparse matrix obtained via Eq.~\eqref{eq:projected}.
The diagonal scaling and iterative solution procedure are presented in
Appendix~\ref{app:fit}.

Hereafter, we denote the construction of sparse matrix entries
from Eq.~\eqref{eq:projected} as \texttt{Truncation}, and denote the approach
in this subsection as \texttt{LeastSquaresFit} or \texttt{LSFit}.

\subsubsection{Approximation quality and iterative convergence}
\label{sec:quality}

The spectrum of matrix $\mathbf{E}_e
=\mathbf{A}_e^{1/2}\mathbf{P}_e\mathbf{A}_e^{1/2}$ is the same
as that of the preconditioned matrix $\mathbf{P}_e\mathbf{A}_e$, and
all eigenvalues equal unity when $\mathbf{P}_e=\mathbf{H}_e$.
The normalised approximation residual from the least-squares problem is
\begin{equation}
 \rho_e
 =
 \frac{1}{\sqrt{N_P}}
 \|\mathbf{I}_{N_P}-\mathbf{E}_e\|_F.
 \label{eq:rho}
\end{equation}
According to Eq.~\eqref{eq:fitspectrum}, $\rho_e$ is the root-mean-square
deviation of the eigenvalues of $\mathbf{E}_e$ from unity. Reducing the fitting
error therefore reduces their aggregate deviation from the
unit spectrum obtained with the exact elemental inverse.

The local spectra can also be related to the global preconditioner.
Let $\mathbf{H}$ and $\mathbf{P}$ denote the
block-diagonal matrices from the exact and approximate
elemental inverses, respectively. Suppose that the eigenvalues
of every $\mathbf{E}_e$ lie in a common interval
$[\alpha,\beta]$, with $0<\alpha\leq\beta$.
In other words, each approximate block satisfies
\begin{equation}
     \alpha\mathbf{H}_e \preceq\mathbf{P}_e \preceq\beta\mathbf{H}_e~.
\end{equation}
Here $\preceq$ denotes the positive-semidefinite ordering.
The global matrices $\mathbf{H}$ and $\mathbf{P}$ satisfy the same inequalities.
For $\mathbf{A}\succ0$, we can get
\begin{equation}
 \kappa_2\!\left(
 \mathbf{A}^{1/2}\mathbf{P}\mathbf{A}^{1/2}
 \right)
 \leq
 \frac{\beta}{\alpha}
 \kappa_2\!\left(
 \mathbf{A}^{1/2}\mathbf{H}\mathbf{A}^{1/2}
 \right)~,
 \label{eq:global-comparison}
\end{equation}
where $\kappa_2(\cdot)$ denotes the condition number in the Euclidean norm.
Thus, if the eigenvalues remain close to unity,
replacing the exact elemental inverses by our approximation 
introduces only a limited deterioration in this condition-number bound.
The least-squares objective does not directly minimise
$\beta/\alpha$, however. It measures the combined error over
all local eigenvalues. A smaller fitting
residual therefore does not necessarily give a smaller global
condition number or fewer solver iterations.
Whether this least-squares approach improves the convergence rate over 
the truncation approach must be assessed in the complete solver and real cases.

Eq.~\eqref{eq:fitting-problem} enforces symmetry, but it does
not ensure the resulting $\mathbf{P}_e$ is positive definite.
Therefore, after constructing $\mathbf{V}_e$ and $\mathbf{P}_e$,
we need to check the positive definiteness, and where necessary, correct it by modifying the diagonal
entries of $\mathbf{V}_e$, which can also increase the fitting residual.
An empirical correction procedure is also briefly discussed in Appendix~\ref{app:fit-spd}.

\section{Numerical results and discussion}
\label{sec:results}
\subsection{Benchmark setup}
\label{sec:case-setup}
\begin{figure}[htbp]
    \centering

    \begin{subfigure}[t]{0.2\textwidth}
        \centering
        \includegraphics[width=\linewidth]{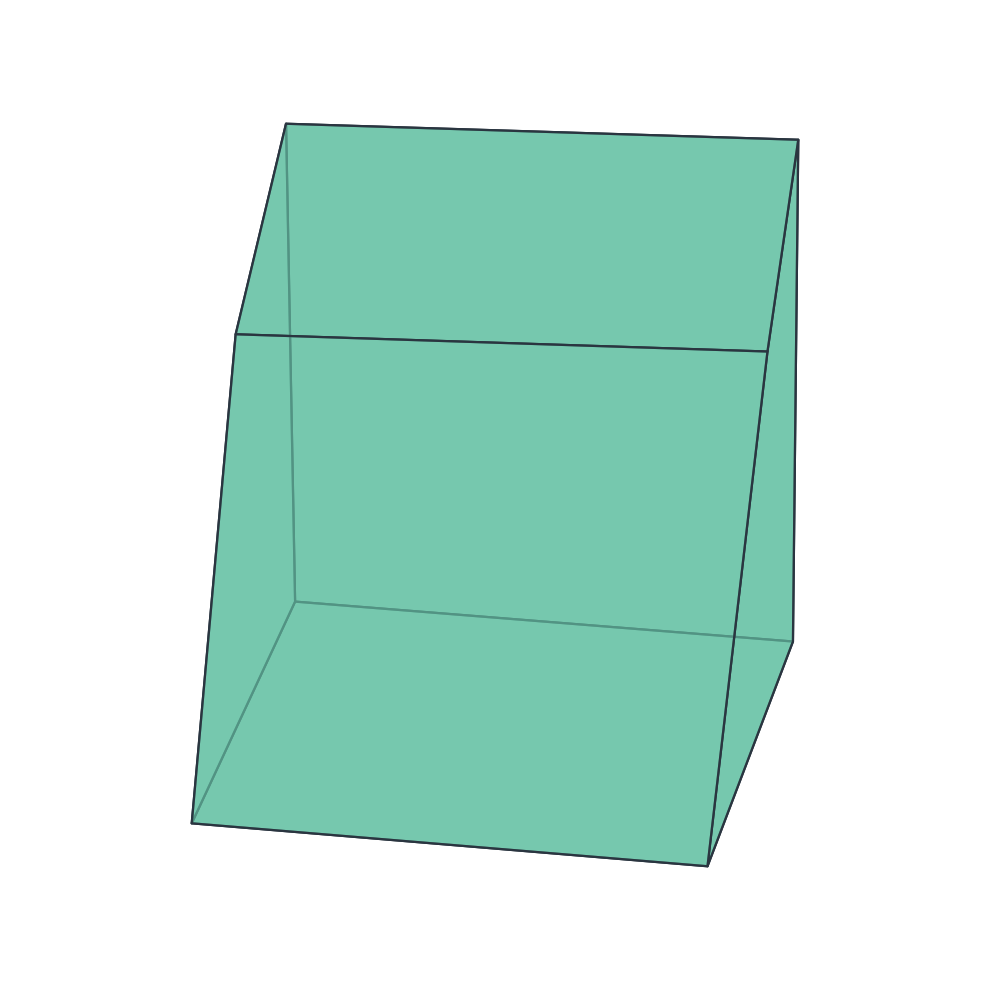}
        \caption{}
        \label{fig:element-hex}
    \end{subfigure}
    \hfill
    \begin{subfigure}[t]{0.2\textwidth}
        \centering
        \includegraphics[width=\linewidth]{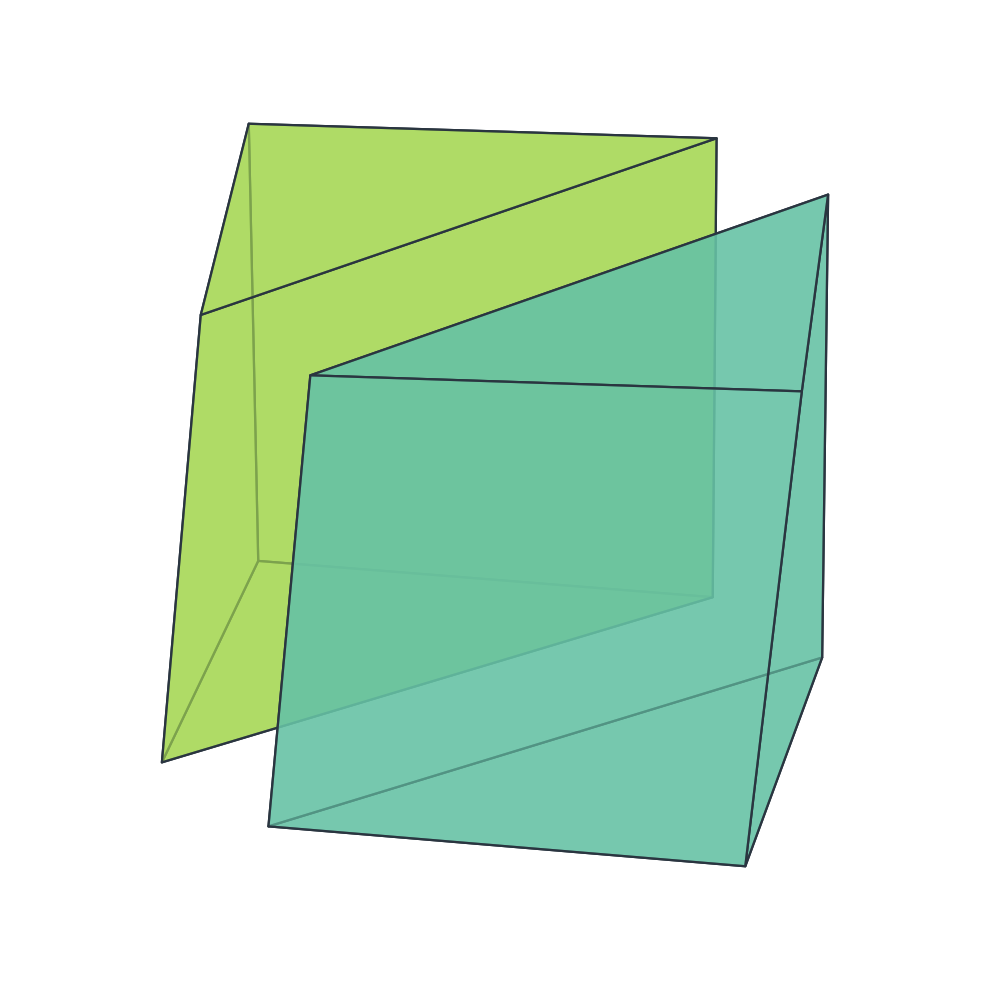}
        \caption{}
        \label{fig:element-prism}
    \end{subfigure}
    \hfill
    \begin{subfigure}[t]{0.2\textwidth}
        \centering
        \includegraphics[width=\linewidth]{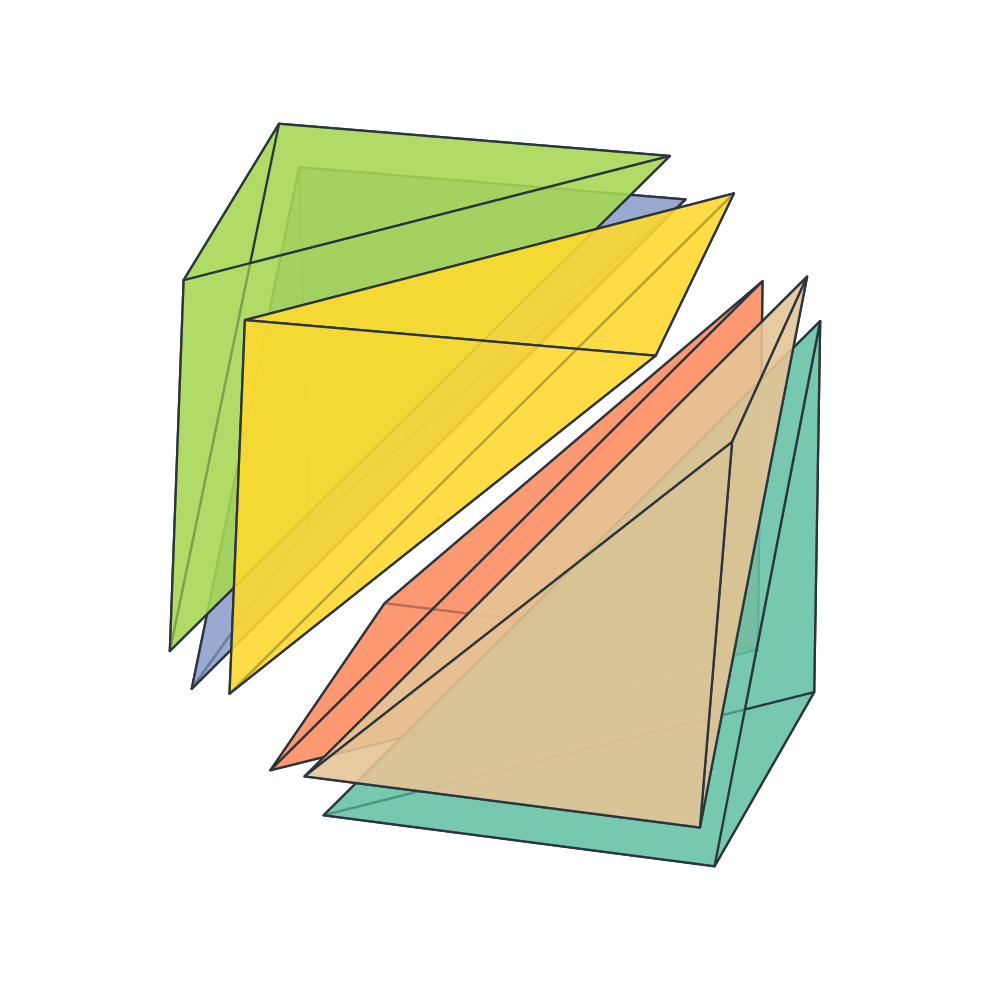}
        \caption{}
        \label{fig:element-tet}
    \end{subfigure}
    \hfill
    \begin{subfigure}[t]{0.2\textwidth}
        \centering
        \includegraphics[width=\linewidth]{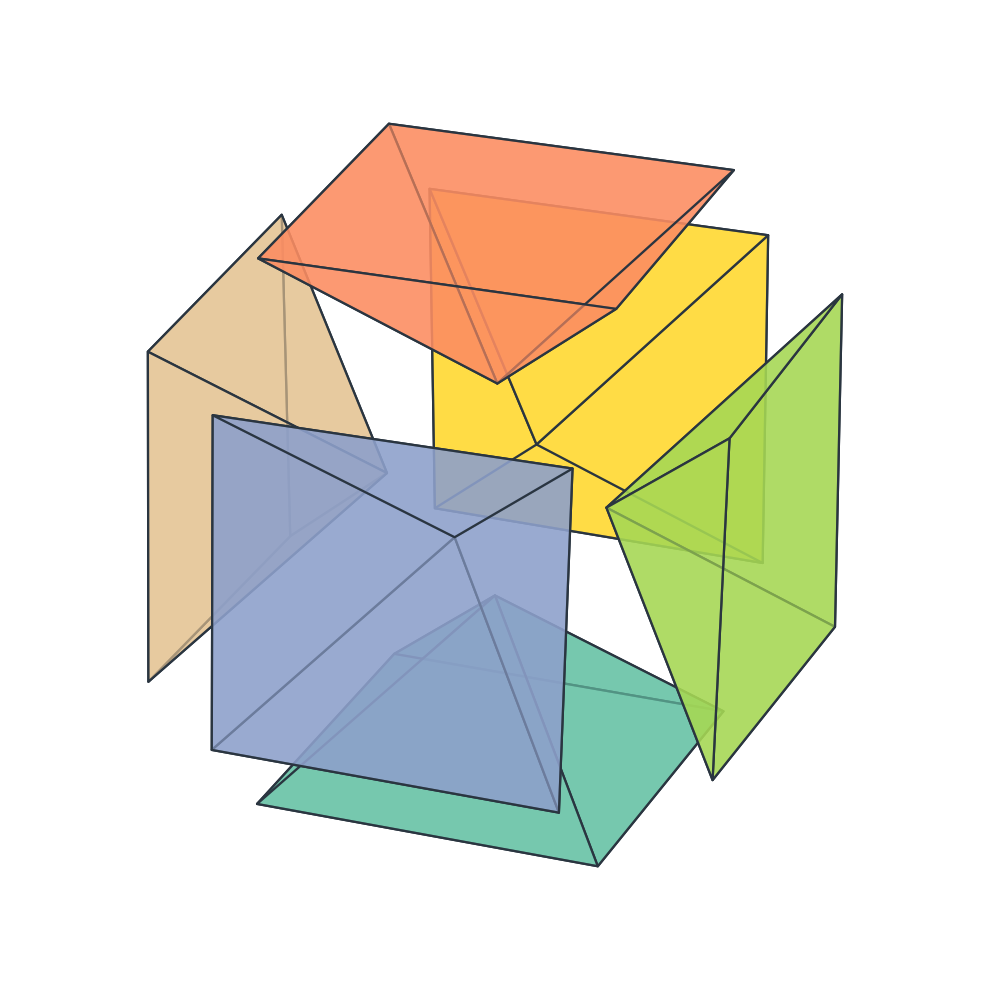}
        \caption{}
        \label{fig:element-pyr}
    \end{subfigure}

    \caption{Different elements in the structured meshes. Each hexahedral element (a) in the base structured
    mesh can be split into 2 prisms (b), 6 tetrahedra (c), or 6 pyramids (d).}
    \label{fig:different-elements}
\end{figure}

We first document the test cases we employ to benchmark the effectiveness and performance of the preconditioners.
We prescribe a manufactured sinusoidal solution for the Poisson problem (see Eq.~\eqref{eq:poisson}) on the cubic domain $\Omega=[-1,1]^3$ as
\begin{equation}
 u_{\text{exact}}=1+\sin(\pi x)\sin(\pi y)\sin(\pi z),\qquad
 f=3\pi^2\sin(\pi x)\sin(\pi y)\sin(\pi z),
\end{equation}
with Dirichlet conditions $g=u_{\text{exact}}$ on all boundaries of the regular cube.
We select this boundary condition because the weak imposition of Dirichlet conditions 
introduces penalty-dependent terms that can significantly affect the conditioning of the discrete system.
The cube domain is discretised using uniform structured base meshes with $N_x^3$ hexahedral cells. 
To form other mesh families, we split each hexahedral cell into either 2 prisms, 6 tetrahedra or 6 pyramids,
as shown in Fig.~\ref{fig:different-elements}. 
The modified modal basis \cite{karniadakis_spectralhp_1999} is used in all these elements. 
We choose three different combinations of mesh size and element order for benchmarks as shown in 
Table~\ref{tab:settings}, reflecting three typical workloads:
low-order fine mesh, medium-order medium-size mesh, and high-order coarse mesh.

\begin{table}[htbp]
\centering
\caption{Mesh families for benchmarks. Each base mesh contains
$N_x^3$ hexahedral cells, and other shapes are subdivided from a hexahedron.
$N_P$ is the number of modal degrees of freedom per
element and $N_{\mathrm{el}}N_P$ is the
total number of DG degrees of freedom (DOFs).}
\label{tab:settings}
\small
\setlength{\tabcolsep}{4pt}
\begin{tabular}{lccrrr}
\toprule
Shape
& $N_{\mathrm{el}}$
& $N_P$
& \shortstack{High-order DOFs\\$p=8,\ N_x=16$}
& \shortstack{Medium-order DOFs\\$p=5,\ N_x=32$}
& \shortstack{Low-order DOFs\\$p=2,\ N_x=64$} \\
\midrule
Hex
& $N_x^3$
& $(p+1)^3$
& $2.986\times10^6$
& $7.078\times10^6$
& $7.078\times10^6$ \\

Prism
& $2N_x^3$
& $\dfrac{(p+1)^2(p+2)}{2}$
& $3.318\times10^6$
& $8.258\times10^6$
& $9.437\times10^6$ \\

Tet
& $6N_x^3$
& $\dfrac{(p+1)(p+2)(p+3)}{6}$
& $4.055\times10^6$
& $1.101\times10^7$
& $1.573\times10^7$ \\

Pyr
& $6N_x^3$
& $\dfrac{(p+1)(p+2)(2p+3)}{6}$
& $7.004\times10^6$
& $1.789\times10^7$
& $2.202\times10^7$ \\
\bottomrule
\end{tabular}
\end{table}


For brevity, we use \texttt{Hex}, \texttt{Prism}, \texttt{Tet} and \texttt{Pyr}
to denote hexahedra, prisms, tetrahedra and pyramids, respectively. 
\texttt{Hex16-p8}, for example, denotes the $16^3$ base mesh at polynomial degree $p=8$.

\begin{figure}[htbp]
    \centering

    \begin{subfigure}[t]{0.3\textwidth}
        \centering
        \includegraphics[width=\linewidth]{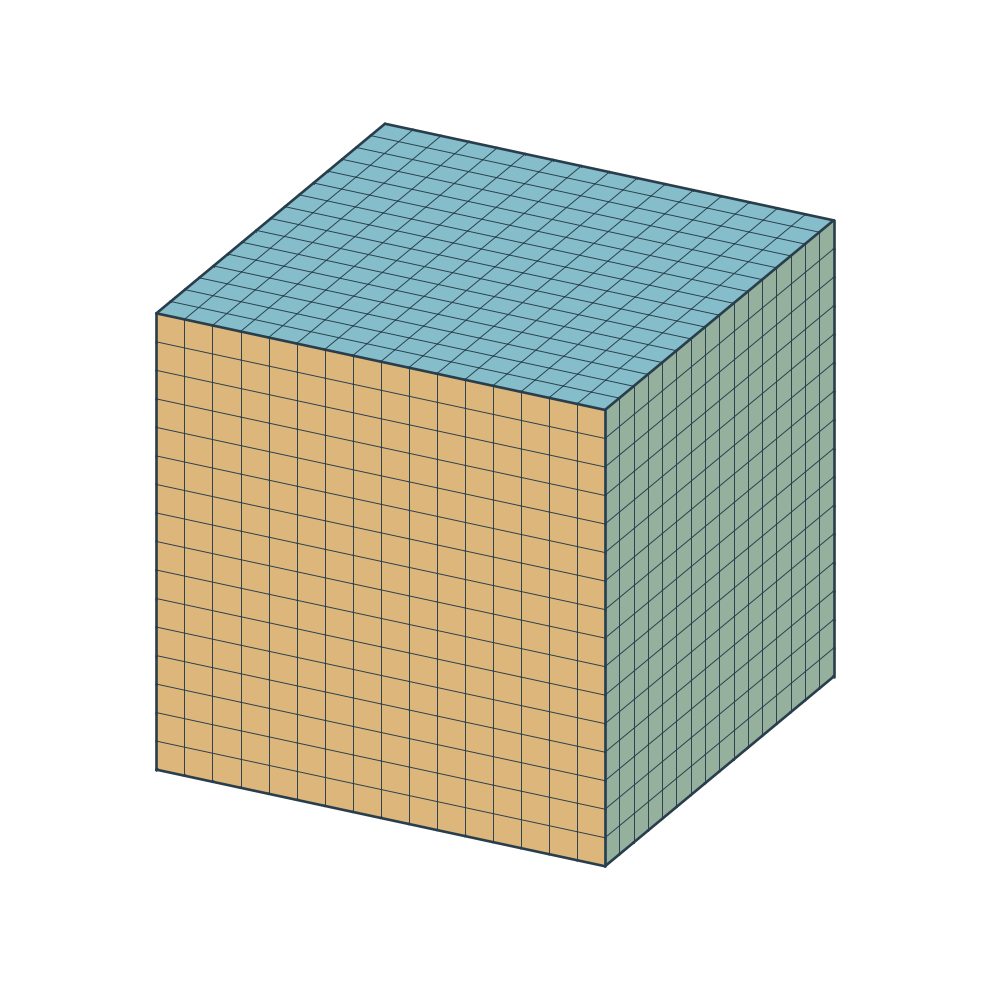}
        \caption{}
        \label{fig:regular-mesh}
    \end{subfigure}
    \begin{subfigure}[t]{0.3\textwidth}
        \centering
        \includegraphics[width=\linewidth]{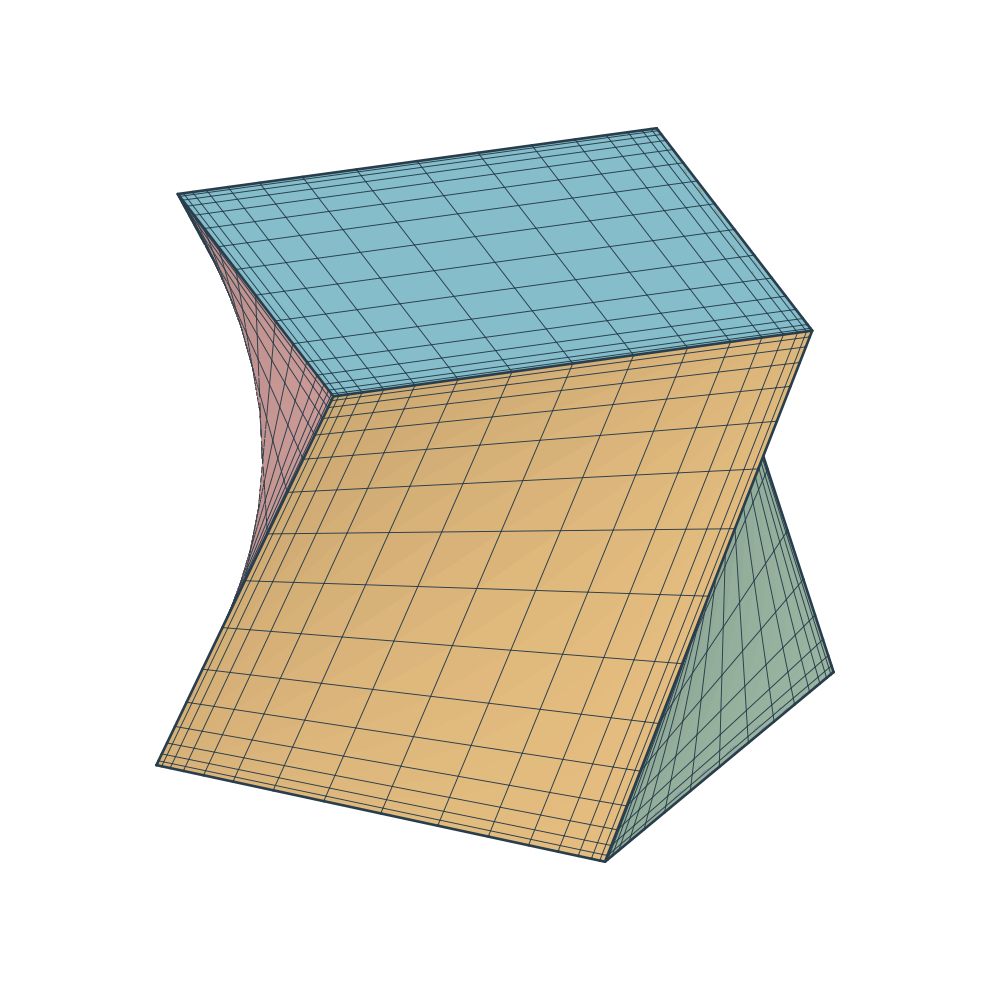}
        \caption{}
        \label{fig:twisted-mesh}
    \end{subfigure}

    \caption{Structured uniform base mesh in a regular cube (a), 
            and the distorted anisotropic structured mesh in a twisted cube (b).}
    \label{fig:regular-twisted-mesh}
\end{figure}

The meshes considered above are uniform, regular, and of good quality. 
In practical applications, however, poor-quality or distorted meshes 
may adversely affect both the accuracy and the convergence of the solver. 
To account for this effect in the benchmarks, 
we twist the cube domain by rotating the top face, and then solve the 
same problem on the
resulting twisted structured meshes of different shapes and mesh sizes.
To introduce additional mesh anisotropy, hyperbolic tangent stretching with the parameter $\beta=2$ is also applied
along each direction. Representative mesh examples are shown in Fig.~\ref{fig:regular-twisted-mesh}.
Hexahedral and prismatic meshes are twisted by 45 degrees, while tetrahedral
and pyramidal meshes are twisted by 30 and 15 degrees, respectively.
Different twisting angles are used to maintain comparable mesh quality across the different element shapes. 
In particular, applying a larger rotation to the pyramidal mesh can produce nearly degenerate elements, 
for which the resulting discrete system becomes difficult to solve effectively.

A preconditioned conjugate gradient solver is used to solve the SIPG Poisson problem, and 
the tolerance of the normalised residual is fixed at $10^{-9}$ in all of the benchmarks.
The polynomial multigrid configurations for both mesh families are summarised in
Table~\ref{tab:mg-settings}. 
The $p$-MG hierarchy for each order is determined by preliminary studies to achieve good overall efficiency. 
For the twisted meshes, smoothing steps and coarse-solver iterations are increased to maintain convergence.
For brievity, we do not optimise the MG configurations individually for each case.

We measure convergence by the number of conjugate gradient iterations executed
to reach the specified tolerance, and measure efficiency
and performance by the elapsed solver time. Setup time is excluded from the timings.

\begin{table}[htbp]
\centering
\caption{Polynomial multigrid settings used in different cases. 
The tuple entries are pre-smoothing steps, post-smoothing steps, 
and coarse-solver iterations, respectively. 
Smoothing applies on each level except the coarsest level.}
\label{tab:mg-settings}
\begin{tabular}{llll}
\toprule
$p$ & $p$-MG hierarchy & \multicolumn{2}{c}{ (pre, post, coarse)} \\
    &                  & Regular mesh & Twisted mesh \\
\midrule
8 & $8\to5\to3\to1$ & $(5,5,5)$     & $(15,15,15)$ \\
5 & $5\to3\to1$     & $(5,5,5)$     & $(15,15,15)$ \\
2 & $2\to1$         & $(5,5,10)$    & $(15,15,15)$ \\
\bottomrule
\end{tabular}
\end{table}

The benchmarks run on a compute node with two AMD EPYC 9554 processors,
providing 128 physical cores in total. 
Each processor has 256~MB of L3 cache and 64~MB of L2 cache,
and uses the AMD Zen 4 architecture with AVX-512 support.
Parallelism is implemented via the Message Passing Interface (MPI). 
Unless stated otherwise, we employ 128 MPI ranks and bind each rank to one physical core.

\subsection{Impact of penalty factor and basis type on convergence}
\label{sec:convergence-factors}

Before comparing preconditioners across element shapes,
we examine tunable factors that may significantly affect the convergence of the SIPG Poisson solver.
This guides the selection of optimal settings for comparing preconditioners in various cases,
and also provide a broader view of convergence issues in the SIPG discretisation.

The first factor is the penalty parameter. As stated in Eq.~\eqref{eq:penalty-parameter} 
in Section~\ref{sec:discretisation}, we employ a shape-dependent penalty formulation
and control the penalty by a constant scaling factor $\sigma_0$.
We vary $\sigma_0$ and repeat the tests
with different case setups, including element shapes, polynomial degrees and mesh sizes.
We only employ the basic block-Jacobi preconditioner.
As shown in Fig.~\ref{fig:penalty-study-regular}, in general, $\sigma_0$ needs to be larger than a critical
value to achieve a converged and accurate solution. 
However, higher $\sigma_0$ can also gradually slow down the convergence. 
Therefore, an optimal penalty scale factor $\sigma_0$ exists, which accounts for both accuracy and convergence, 
and is roughly independent of the mesh size and orders. For all the shapes,
the optimal $\sigma_0$ is roughly between 0.2 and 0.25, while an exception is the $p=1$ tetrahedral case,
which requires approximately 0.32 for stable convergence in these tests.
Fig.~\ref{fig:penalty-study-twisted} shows the results obtained on twisted meshes,
and we find that the optimal $\sigma_0$ does not change much in these cases.
In our benchmarks, we set $\sigma_0=0.32$ in all the regular mesh cases and $\sigma_0=0.4$ in 
all the distorted mesh cases, leaving sufficient margin.

\begin{figure}[ht]
    \centering
    \includegraphics[width=0.75\linewidth]{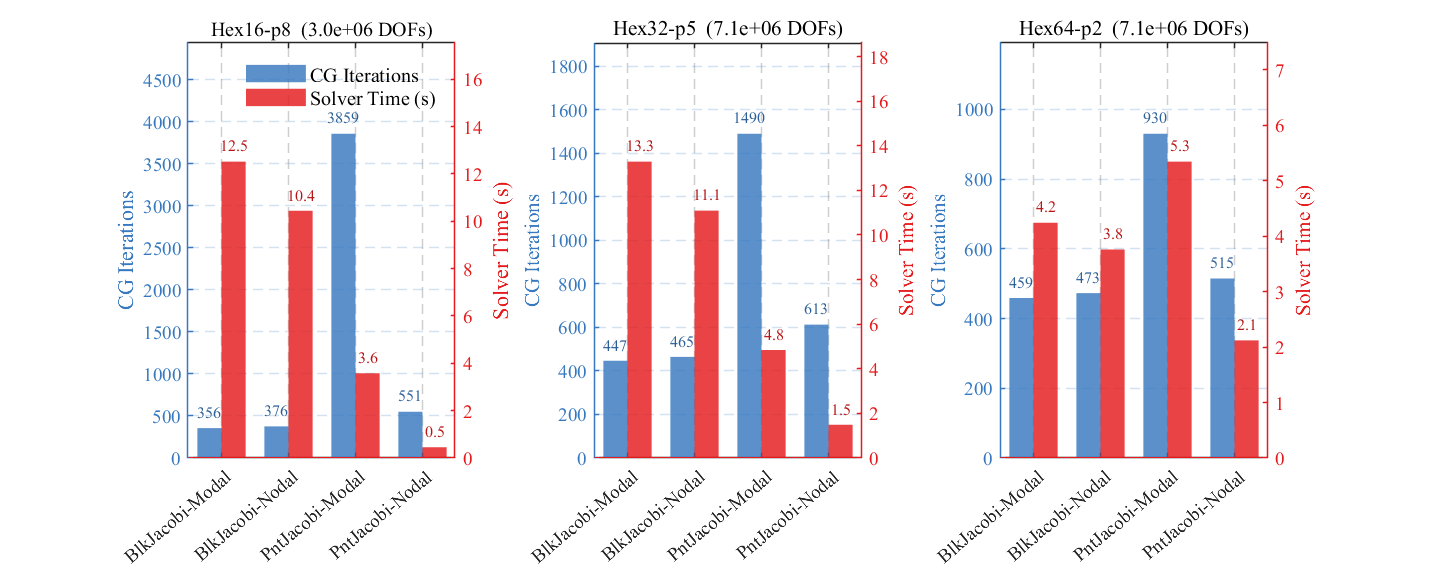}
    \caption{Conjugate gradient iterations and solver time required to achieve a tolerance of 
    $10^{-9}$ for the regular hexahedral cases 
    with different basis types (modal and nodal) and preconditioners (point Jacobi and block Jacobi).}
    \label{fig:precon-comparison-nodal-modal-hex}
\end{figure}

The second factor is the choice of spectral element basis.
Fig.~\ref{fig:precon-comparison-nodal-modal-hex} compares the convergence results
from two common choices of basis in hexahedral meshes.
One is the popular Lagrange basis using tensor-product GLL points, usually called the nodal basis.
The other is the modal basis used in Nektar++.
When we use a point-Jacobi preconditioner, 
the nodal case converges much faster than the modal case, especially for high orders.
If we change to a block-Jacobi preconditioner,
the iteration counts for the nodal and modal bases are similar,
but the lower iteration count does not always reduce solver time 
because the cost per call of block Jacobi is much higher than that of point Jacobi.
The tensor-product GLL nodal basis gives better conditioning with point Jacobi, 
but this basis is only available for quadrilaterals and hexahedra.
As stated earlier in Section~\ref{sec:construction},
these findings motivated us to use a modal-to-nodal transform
when constructing the approximate inverse preconditioner for the modal basis.
In the following, we will focus on the convergence of the SIPG Poisson solver with a modal basis 
and show how we can improve the convergence and performance by constructing
a matrix-free approximate elemental inverse preconditioner. The aim is to obtain
iteration counts between those of point Jacobi and block Jacobi,
while achieving a shorter solver time than block Jacobi.

\subsection{Selecting the fastest approximate elemental inverse preconditioner}
\label{sec:fastest-approx-inverse}
We perform a series of tests for each configuration of the approximate elemental inverse 
to find the best combination.
The configurable components are:
\begin{itemize}
    \item Choice of projection space (see Section~\ref{sec:projection});
    \item Choice of sparse pattern (see Section~\ref{sec:sparse-pattern});
    \item Implementation of the inverse transform (see Section~\ref{sec:inverse-sum-factorisation});
    \item Construction of sparse matrix entries (see Section~\ref{sec:least-squares-fit});
\end{itemize}
For readability, we present selected results when comparing specific
configurations. The complete results are given in
Figs.~\ref{fig:sipg-preconditioners-hex}--\ref{fig:sipg-preconditioners-pyr}
in Appendix~\ref{app:results}.

Table~\ref{tab:different-projection} compares the different choices of projection spaces.
The tensorial GLL nodal space gives the lowest iteration count in all cases shown.
The original quadrature space has similar convergence, but it is bigger
than the GLL nodal space and incurs additional computational cost.
The non-tensorial nodal space gives poorer convergence in these tests,
despite its symmetric node distribution and smaller size. Therefore, in the following discussions,
the comparisons in this section use \texttt{ModalToGLL} projection.

\begin{table}[htb]
    \caption{Conjugate gradient (CG) iteration counts for different projection spaces with the approximate elemental
    inverse preconditioner. All cases use regular meshes, the \texttt{Band7}
    sparse pattern, the dense pseudoinverse transform (\texttt{FullPseudoInv}), and sparse
    matrices obtained by \texttt{Truncation}.}
    \centering
    \begin{tabular}{lcccc}
    \toprule
        \textbf{Mesh} & $p$ & \multicolumn{3}{c}{\textbf{CG iterations}} \\
    \cmidrule(l){3-5}
                     &     & \texttt{ModalToGLL} & \texttt{ModalToPhys} & \texttt{ToNonTensorNodal} \\
    \midrule
        Hex16   & 8 &  568 &  761 &     /       \\
        Hex32   & 5 &  596 &  794 &     /       \\
    \cmidrule{1-5}
        Prism16 & 8 &  770 & 1049 & 1614 \\
        Prism32 & 5 &  778 & 1065 & 937 \\
    \cmidrule{1-5}
        Tet16   & 8 & 1459 & 1473 & 4076 \\
        Tet32   & 5 & 1410 & 1471 & 1583 \\
    \cmidrule{1-5}
        Pyr16   & 8 & 1076 & 1437 &     /       \\
        Pyr32   & 5 & 1209 & 1481 &     /       \\
    \bottomrule
    \end{tabular}
    \label{tab:different-projection}
\end{table}

\begin{table}[htb]
    \caption{The influence of different inverse transforms on the convergence of the conjugate gradient
    solver with the approximate elemental inverse preconditioner. 
    All cases use regular meshes, the \texttt{Band7} sparse pattern, 
    and use sparse matrices obtained by \texttt{Truncation}.
    Each entry gives the number of CG iterations followed by the solver time in seconds.
    The fastest solver time is underlined.}
    \centering
    \begin{tabular}{llccc}
    \toprule
        \textbf{Mesh} & $p$ & \multicolumn{3}{c}{\textbf{CG iterations (solver time, s)}} \\
    \cmidrule(lr){3-5}
                      &     & \texttt{FullPseudoInv} & \texttt{SumFacPseudoInv} & \texttt{SumFacLeftInv} \\
    \midrule
        Hex16         & 8   & 568 (15) & 568 \underline{(0.67)}  &     /                 \\
        Hex32         & 5   & 596 (11) & 596 \underline{(2.5)}  &     /                 \\
    \cmidrule{1-5}
        Prism16       & 8   & 770 (22) & 770 \underline{(2.2)}  & 868 (2.4)             \\
        Prism32       & 5   & 778 (18) & 778 \underline{(5.5)}  & 894 (6.4)             \\
    \cmidrule{1-5}
        Tet16         & 8   & 1459 (51) &       /           & 2471 \underline{(19)} \\
        Tet32         & 5   & 1410 (58) &       /          & 2035 \underline{(43)} \\
    \cmidrule{1-5}
        Pyr16         & 8   & 1076 (64) &       /           & 2287 \underline{(21)} \\
        Pyr32         & 5   & 1209 (70) &       /          & 2255 \underline{(56)} \\
    \bottomrule
    \end{tabular}
    \label{tab:different-inverse}
\end{table}

The choice of inverse transform strongly depends on element shapes.
For hexahedral elements, the \texttt{SumFacPseudoInv} is mathematically
equivalent to \texttt{FullPseudoInv} but has significantly lower computational cost.
Therefore, we should always use \texttt{SumFacPseudoInv} in hexahedral elements.
For tetrahedral and pyramidal elements, the \texttt{SumFacLeftInv}
is not equivalent to \texttt{FullPseudoInv}. We need to examine if
the left inverse harms the convergence and whether
the lower application cost compensates for the additional iterations.
For prismatic elements, both \texttt{SumFacPseudoInv} and
\texttt{SumFacLeftInv} exist; we also need to compare which version
of sum-factorisation gives the best outcome.
Table~\ref{tab:different-inverse} displays representative results to show the
influence of different inverse transforms. 
In general, \texttt{SumFacLeftInv} slows convergence, compared with
a pseudoinverse. For the fixed \texttt{Band7}/\texttt{Truncation} configurations in this table,
\texttt{SumFacPseudoInv} is preferred on prisms.
For tetrahedra and pyramids, \texttt{SumFacLeftInv} is the only available
sum-factorised option. Although it requires more iterations, the conjugate gradient solver time
is still less than that of \texttt{FullPseudoInv}, thanks to the sum-factorisation.
Higher-order cases can cause a more significant performance gap, as shown in the full results
in Appendix~\ref{app:results}.

\begin{table}[htb]
    \caption{The influence of least-squares fitting on the convergence of the conjugate gradient
    solver with the approximate elemental inverse preconditioner.
    All cases use regular meshes, project to GLL nodal space (\texttt{ToGLL}),
    and use the \texttt{Band7} sparse pattern.}
    \centering
    \begin{tabular}{lllcc}
    \toprule
        \textbf{Mesh} & $p$  & \textbf{Inverse transform} & \multicolumn{2}{c}{\textbf{CG iterations}}  \\
    \cmidrule{4-5}
                      &      &                            &  \texttt{Truncation} & \texttt{LeastSquaresFit} \\
    \midrule
        Hex16   &  8         & \texttt{SumFacPseudoInv}  &         568          &       519    \\
        Hex32   &  5         &                           &         596          &       600      \\
    \cmidrule{1-5}
        Prism16 &  8         & \texttt{SumFacPseudoInv}  &         770         &       741     \\
        Prism32 &  5         &                           &         778         &       806     \\
    \cmidrule{1-5}        
        Tet16   &  8         & \texttt{SumFacLeftInv}    &          2471        & 1393     \\
        Tet32   &  5         &                           &          2035        &  1226    \\
    \cmidrule{1-5}         
        Pyr16   &  8         & \texttt{SumFacLeftInv}    &          2287        & 1057    \\
        Pyr32   &  5         &                           &          2255        & 1214    \\
    \bottomrule
    \end{tabular}
    \label{tab:least-squares-convergence}
\end{table}

\texttt{LeastSquaresFit} can be applied after \texttt{Truncation}
to reduce the local fitting residual.
Table~\ref{tab:least-squares-convergence} presents the convergence rates before and after applying the least-squares fitting.
For medium-to-high-order tetrahedral and pyramidal cases, \texttt{LeastSquaresFit} is particularly
important to recover the convergence rate with \texttt{SumFacLeftInv} enabled.
For hexahedral and prismatic meshes, however, \texttt{LeastSquaresFit} does not consistently improve convergence
when using \texttt{SumFacPseudoInv}. In other words, \texttt{SumFacPseudoInv} together with the
\texttt{Truncation} in Eq.~\eqref{eq:projected} inherently produces a good estimate of the sparse entries
so that further optimisation makes little difference. It should be noted that
applying additional \texttt{LeastSquaresFit} may significantly increase the solver setup time,
especially in high-order cases. Table~\ref{tab:least-squares-setup-time} shows the preconditioner setup
time in tetrahedral cases. Solving an additional least-squares problem for each element can increase
the setup time by more than a factor of ten. This overhead can be mitigated in at least two ways. First, the initialisation
is performed once for a fixed discretisation, and the cost can be amortised over thousands of repeated solves.
Second, elements can have similar or identical local operators; thus, we do not have to re-solve the least-squares
problem for identical elements, or solve each problem from scratch.

\begin{table}[htb]
    \caption{The influence of least-squares fitting on the setup time of the preconditioner. 
        All cases project to GLL nodal space (\texttt{ToGLL}) and use the \texttt{Band7} sparse pattern. 
        Each elemental block inverse and least-squares fitting
        is computed independently from scratch.}
    \centering
    \begin{tabular}{llcc}
    \toprule
        \textbf{Mesh} & $p$ & \multicolumn{2}{c}{\textbf{Preconditioner setup time (s)}} \\
    \cmidrule{3-4}
                      &     & \texttt{Truncation} & \texttt{LeastSquaresFit} \\
    \midrule
        Tet16         &  8  & 8.2  & 126 \\
        Tet32         &  5  & 2.0  & 39  \\
        Tet64         &  2  & 0.13 & 2.9  \\
    \bottomrule
    \end{tabular}
    \label{tab:least-squares-setup-time}
\end{table}

The sparse pattern is generally a secondary choice and less influential on the 
performance and convergence once the projection space and the inverse transform are selected. 
As documented in the appendix and in Table~\ref{tab:different-sparse-pattern}, 
the influence of sparse pattern varies according to the shapes and orders. 
For hexahedra, the simplest diagonal matrix (\texttt{Band1}) gives
the shortest solver time in Table~\ref{tab:different-sparse-pattern}.
Denser patterns can reduce iteration counts, but their higher application
cost outweighs this improvement in the cases shown. However, for other collapsed-coordinate elements,
retaining neighbouring-point couplings can be valuable.
For prisms, the \texttt{NearestPnt} pattern gives the shortest solver time for the high-order cases shown,
while the \texttt{Band1} sparse pattern gives the fastest solution for medium- or low-order cases.
For tetrahedra and pyramids, the conclusion is similar. 
The \texttt{NearestPnt} pattern gives the shortest solver time at $p=8$ and $p=5$,
while \texttt{Band1} gives the shortest solver time at $p=2$.
Table~\ref{tab:different-sparse-pattern} also marks the fastest configuration
of the approximate elemental inverse preconditioner for the different workloads.
If we choose not to apply least-squares optimisation, then the best configuration may slightly
change, which can be read from the full results in the appendix.

\begin{table}[htb]
    \centering
    \caption{The influence of sparse pattern on the convergence of the conjugate gradient solver
        with the approximate elemental inverse preconditioner. 
        All cases use regular meshes, project to GLL nodal space (\texttt{ToGLL}) 
        and obtain sparse matrices by \texttt{LeastSquaresFit}.
        Each entry gives the number of outer CG iterations and solver time in seconds (in parentheses).
        The fastest solver time in each row is underlined.}
    \begin{tabular}{lllcccc}
    \toprule
        \textbf{Mesh} & $p$ &  \textbf{Inverse transform}   & \multicolumn{4}{c}{ \textbf{CG iterations (solver time, s)}  }                        \\
    \cmidrule(lr){4-7}           
                      &     &                               & \texttt{Band1}  &  \texttt{Band7}    & \texttt{NearestPnt}  & \texttt{FirstLayer}  \\
    \midrule                                
        Hex16         &  8  &                               &    519 \underline{(0.54)}  &     519 (0.61)     &      456 (0.57)     &      420 (0.77)        \\
        Hex32         &  5  &   \texttt{SumFacPseudoInv}    &    590 \underline{(2.0)}   &     600 (2.5)      &      675 (2.7)      &      533 (2.9)         \\
        Hex64         &  2  &                               &    482 \underline{(2.8)}   &     582 (3.8)      &      662 (4.2)      &      612 (4.4)         \\
    \cmidrule{1-7}                               
        Prism16       &  8  &                               &    770 (1.8)            &     741 (2.0)      &      609 \underline{(1.7)}     &      606 (2.7)         \\
        Prism32       &  5  &    \texttt{SumFacPseudoInv}   &    835 \underline{(4.8)}   &     806 (5.7)      &      734 (5.1)      &      722 (6.9)         \\
        Prism64       &  2  &                               &    705 \underline{(7.3)}   &     799 (9.4)      &      904 (10)       &      881 (11)          \\
    \cmidrule{1-7}                                   
        Tet16         &  8  &                               &   2098 (13)             &    1393 (11)       &      837 \underline{(6.5)}     &      811 (9.8)         \\
        Tet32         &  5  &   \texttt{SumFacLeftInv}      &   1586 (28)             &    1226 (26)       &     1095 \underline{(22)}      &      996 (28)          \\
        Tet64         &  2  &                               &   1243 \underline{(45)}    &    1305 (52)       &     1320 (51)       &     1202 (53)          \\
    \cmidrule{1-7}                                    
        Pyr16         &  8  &                               &   1461 (11)             &    1057 (9.7)      &      801 \underline{(7.3)}     &      747 (9.5)         \\
        Pyr32         &  5  &   \texttt{SumFacLeftInv}      &   1396 (30)             &    1214 (30)       &     1112 \underline{(27)}      &     1124 (35)          \\
        Pyr64         &  2  &                               &   1103 \underline{(49)}   &    1133 (54)      &     1336 (62)       &     1207 (61)          \\
    \bottomrule
    \end{tabular}
    \label{tab:different-sparse-pattern}
\end{table}

\begin{figure}[ht]
    \centering
    \includegraphics[width=0.97\linewidth]{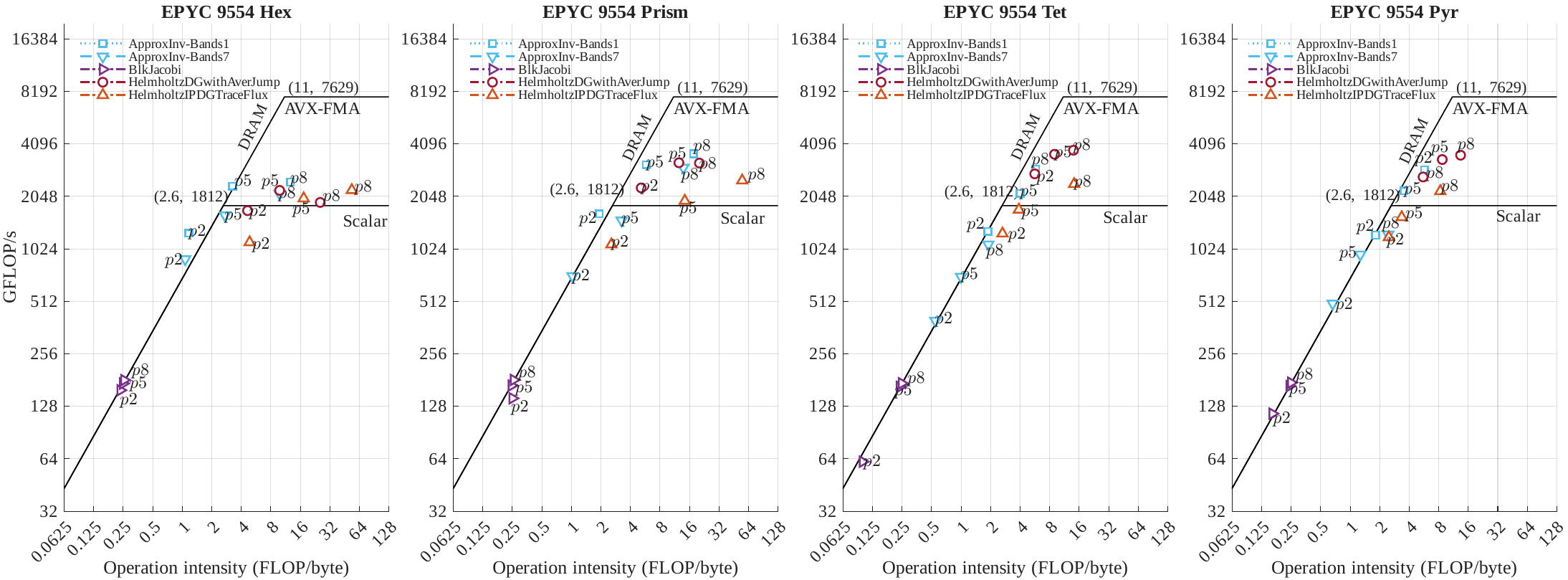}
    \caption{Roofline results of the preconditioners on the dual-socket EPYC 9554 compute node. Each category has three points
        representing three workloads in Table~\ref{tab:settings}.
        \texttt{HelmholtzDGwithAverJump} and \texttt{HelmholtzIPDGTraceFlux} are the major
         matrix-free operators for evaluating the SIPG formulation.}
    \label{fig:sipg-precon-roofline}
\end{figure}

To validate the efficiency of the matrix-free preconditioners, 
we collected profiling results using \texttt{likwid}~\cite{thomas_gruber_2024_likwid},
and show the roofline plots for different preconditioners on different workloads 
in Fig.~\ref{fig:sipg-precon-roofline}.
The two major customisable components in a conjugate gradient solver are 
the left-hand side evaluation ($\mathbf{A}\mathbf{x}$) and the preconditioning.
Our preceding work \cite{xia2026} tackles the first component, and we have identified
and constructed two fused matrix-free operators named \texttt{HelmholtzDGwithAverJump}
and \texttt{HelmholtzIPDGTraceFlux} for evaluating the left-hand side (LHS). In the current work,
we tackle the second component, leading to the matrix-free approximate inverse preconditioner. 
As can be seen in Fig.~\ref{fig:sipg-precon-roofline}, the approximate inverse preconditioner,
which projects to the GLL nodal space and uses a sum-factorisation inverse,
successfully exploits the SIMD capability on the tested platform, and achieves performance
similar to that of the matrix-free operators for LHS evaluation.
Higher polynomial degrees lead to higher operation intensity and higher FLOPS (floating-point operations per second),
making both the LHS and preconditioner more advantageous.
In contrast, the block-Jacobi preconditioner has a much lower operation intensity around $1/4$ FLOP/byte, which
is expected, as GEMV executes 2 floating-point operations for each matrix entry loaded (8 bytes).

\subsection{Applying the approximate elemental inverse to multigrid smoothing}
\label{sec:mg-results}

The preceding tests show that the approximate elemental inverse can reduce the solver time
compared with the block-Jacobi or point-Jacobi preconditioner.
Now we apply it to the Chebyshev smoothers in our multigrid preconditioner to see if
better performance can be achieved. Figs.~\ref{fig:sipg-continuous-hex}--\ref{fig:sipg-continuous-pyr} in Appendix~\ref{app:results}
present the convergence and performance of the SIPG Poisson solver with multigrid preconditioners.
The multigrid procedure has been documented in Section~\ref{sec:mg}, and the hierarchy,
smoothing steps, and coarse-solver iterations follow Table~\ref{tab:mg-settings}.
For each workload and mesh family, we select the approximate inverse with the shortest
multigrid solver time among the configurations shown in these figures.
Roughly speaking, the best configuration is the same as in the standalone preconditioner.
Table~\ref{tab:mg-smoother-comparison} summarises the results and compares them with
point-Jacobi and block-Jacobi smoothing under the same multigrid settings.

We first examine the results on regular meshes. Compared with the standalone
preconditioner, multigrid substantially reduces the number of conjugate gradient iterations,
but the cost of each iteration is high due to multi-step smoothing and additional overheads.
Whether the solver time can be reduced effectively depends on the multigrid configurations,
and the properties of the system. The choice of smoother makes a large difference. 
Point Jacobi is not an effective smoother, particularly in the high-order cases.
The resulting SIPG solver is even slower than using a standalone point-Jacobi preconditioner.
The block-Jacobi and approximate inverse smoothers achieve a similar number of outer iterations,
but the approximate inverse is significantly cheaper to apply, and therefore can give shorter solver time. 
The performance advantage of the approximate inverse over block Jacobi
is most significant in high-order hexahedra and prisms, 
with an order of magnitude difference in the $p=8$ hexahedral case.
For tetrahedral and pyramidal meshes, the advantage of the approximate inverse over block Jacobi is less significant.
For the low-order cases ($p=2$), the advantage over block Jacobi is modest.

\begin{table}[htbp]
    \centering
    \caption{Convergence and performance of the SIPG Poisson solver with different multigrid smoothers.
        Each entry gives the number of outer CG iterations and solver time in seconds (in parentheses).
        The approximate inverse uses the fastest configuration reported in the
        corresponding panel of Appendix~\ref{app:results}.
        The shortest solver time for each workload is underlined.}
    \label{tab:mg-smoother-comparison}
    \small
    \setlength{\tabcolsep}{4pt}
    \begin{tabular}{llcccccc}
    \toprule
        \textbf{Mesh} & $p$ & \multicolumn{3}{c}{\textbf{Regular mesh}} & \multicolumn{3}{c}{\textbf{Twisted mesh}} \\
    \cmidrule(lr){3-5}\cmidrule(lr){6-8}
                      &     & \texttt{PointJacobi} & \texttt{BlockJacobi} & \texttt{ApproxInverse} & \texttt{PointJacobi} & \texttt{BlockJacobi} & \texttt{ApproxInverse} \\
    \midrule
        Hex16   & 8 & 478 (12.0) & 7 (2.19) & 8 \underline{(0.21)} & 292 (27.4) & 7 (6.81) & 8 \underline{(0.76)} \\
        Hex32   & 5 & 177 (13.0) & 7 (2.05) & 8 \underline{(0.61)} & 164 (59.3) & 6 (6.49) & 8 \underline{(2.96)} \\
        Hex64   & 2 & 106 (26.7) & 7 (2.00) & 7 \underline{(1.79)} &  74 (50.7) & 8 (6.35) & 8 \underline{(5.50)} \\
    \cmidrule{1-8}
        Prism16 & 8 & 318 (15.3) & 8 (2.08) & 9 \underline{(0.53)} & 310 (60.7) & 8 (6.45) & 9 \underline{(2.06)} \\
        Prism32 & 5 & 164 (19.5) & 8 (2.32) & 8 \underline{(1.08)} & 209 (131)  & 8 (8.82) & 9 \underline{(6.08)} \\
        Prism64 & 2 &  70 (23.1) & 8 (2.91) & 8 \underline{(2.73)} &  85 (90.2) & 11 (12.6) & 11 \underline{(11.8)} \\
    \cmidrule{1-8}
        Tet16   & 8 & 274 (36.4) & 9 (2.11) & 13 \underline{(2.02)} & 173 (50.1) & 8 \underline{(4.78)} & 37 (13) \\
        Tet32   & 5 & 166 (62.2) & 8 \underline{(3.91)} & 11 (4.67) & 141 (118) & 9 (10.8) & 10 \underline{(9.72)} \\
        Tet64   & 2 &  62 (50.3) & 11 (9.49) & 10 \underline{(8.47)} & 52 (77.5) & 9 \underline{(15.0)} & 11 (17.6) \\
    \cmidrule{1-8}
        Pyr16   & 8 & 426 (66.6) & 9 (4.32) & 10 \underline{(2.32)} & 233 (165) & 10 (17.5) & 10 \underline{(9.62)} \\
        Pyr32   & 5 & 235 (109)  & 8 (5.82) & 9 \underline{(4.64)} & 165 (344) & 10 (28.7) & 10 \underline{(24.9)} \\
        Pyr64   & 2 & 117 (125)  & 9 (10.4) & 9 \underline{(9.95)} &  62 (236) & 10 (40.1) & 10 \underline{(39.0)} \\
    \bottomrule
    \end{tabular}
\end{table}

Table~\ref{tab:mg-smoother-comparison} also shows how the multigrid performance changes
on the twisted meshes. For hexahedra and prisms, the outer iteration count is similar to that on regular meshes,
but solver time increases partly because the twisted cases use 15 smoothing steps
rather than 5,
as shown in Table~\ref{tab:mg-settings}.
The approximate inverse still gives convergence close to block Jacobi and retains a shorter solver time,
but the performance gain becomes smaller. 
For tetrahedra and pyramids, the iteration counts and solver times are
more sensitive to distortion, with the changes depending on degree and smoother.
An outlier is the twisted $p=8$ tetrahedral case, where
the approximate inverse requires substantially more outer iterations than block Jacobi.
Denser patterns also remain slower than exact block Jacobi in this case.
The approximate inverse is not guaranteed to be the best choice in
collapsed-coordinate elements and distorted meshes, 
and we need to check before applying it in realistic scenarios.

\subsection{Comparison with the continuous Galerkin solver}
\label{sec:sipg-continuous-comparison}

Figs.~\ref{fig:sipg-continuous-hex}--\ref{fig:sipg-continuous-pyr}
also include the benchmark results from the continuous Galerkin solver in Nektar++, 
which either solves a full system or a statically condensed system. 
The full system retains both element-interior and element-boundary DOFs, and solves them together. 
Static condensation \cite{kirby2012CGHDG} eliminates the interior DOFs by substitution locally, 
leaving a smaller global Schur-complement system for the element-boundary DOFs, 
which can be solved more easily and cheaply.
The interior solution is then recovered by local back substitution.
In addition to the point-Jacobi and block preconditioners,
a low-energy block preconditioner \cite{sherwin2001} is also available for 
the statically condensed continuous Galerkin system in Nektar++,
which is more effective than the block preconditioner.
Table~\ref{tab:sipg-continuous-comparison} selects
the shortest solver time among these continuous Galerkin configurations and compares
it with the fastest SIPG choice in Table~\ref{tab:mg-smoother-comparison}.

In the comparison figures, \texttt{Galerkin-Full} and
\texttt{Galerkin-StaticCond} denote the full and statically condensed
continuous Galerkin systems. \texttt{Block} and \texttt{LowEngBlk} denote
the block and low-energy block preconditioners, respectively.

SIPG with $p$-MG and the approximate inverse smoother is faster than
continuous Galerkin for the medium-to-high-order regular hexahedral, prismatic
and pyramidal meshes. Continuous Galerkin remains faster for all tetrahedral
and low-order workloads shown.
With a low-energy block preconditioner, the convergence of the continuous Galerkin solver 
is not much affected by the twisted meshes, and thus it is faster than SIPG in most cases. 
The continuous Galerkin solver enforces continuity directly without a significant penalty term
and has a much smaller system size than DG
due to continuity constraints and static condensation.
Therefore, it is not surprising to see continuous Galerkin faster than SIPG.
\citet{kronbichler2018CGvsDG} reached a similar conclusion for smooth solutions on
quadrilateral and hexahedral meshes: with state-of-the-art multigrid solvers,
high-order continuous elements gave the best time to solution, closely followed
by matrix-free SIPG. Unlike that study, the continuous Galerkin solver here
uses static condensation with a low-energy block preconditioner rather than multigrid,
and the present comparison also covers prismatic, tetrahedral and pyramidal elements.

\begin{table}[htbp]
    \centering
    \caption{Comparison of the fastest SIPG solver and fastest continuous Galerkin
        solver for different workloads, including regular and twisted meshes. 
        The entries give solver times in seconds, excluding setup time.
        The outer CG tolerance is fixed at $10^{-9}$ for all cases.}
    \label{tab:sipg-continuous-comparison}
    \small
    \begin{tabular}{llcccc}
    \toprule
        \textbf{Mesh} & $p$ & \multicolumn{2}{c}{\textbf{Regular mesh}} & \multicolumn{2}{c}{\textbf{Twisted mesh}} \\
    \cmidrule(lr){3-4}\cmidrule(lr){5-6}
                      &     & SIPG & Galerkin & SIPG & Galerkin \\
    \midrule
        Hex16   & 8 & 0.21 & 0.72 & 0.76 & 0.89 \\
        Hex32   & 5 & 0.61 & 1.4 & 3.0 & 2.2 \\
        Hex64   & 2 & 1.8 & 0.59 & 5.5 & 1.1 \\
    \cmidrule{1-6}
        Prism16 & 8 & 0.53 & 1.3 & 2.1 & 1.3 \\
        Prism32 & 5 & 1.1 & 2.9 & 6.1 & 3.2 \\
        Prism64 & 2 & 2.7 & 0.87 & 12 & 1.5 \\
    \cmidrule{1-6}
        Tet16   & 8 & 2.0 & 1.3 & 4.8 & 1.3 \\
        Tet32   & 5 & 3.9 & 2.6 & 9.7 & 3.2 \\
        Tet64   & 2 & 8.5 & 1.1 & 15 & 1.0 \\
    \cmidrule{1-6}
        Pyr16   & 8 & 2.3 & 4.0 & 9.6 & 4.3 \\
        Pyr32   & 5 & 4.6 & 5.4 & 25 & 9.7 \\
        Pyr64   & 2 & 10 & 2.6 & 39 & 3.2 \\
    \bottomrule
    \end{tabular}
\end{table}

\FloatBarrier
\section{Conclusions}
\label{sec:conclusions}

We have presented a matrix-free approximate elemental inverse for
preconditioning discontinuous Galerkin spectral element discretisations.
The construction combines inverse modal-to-nodal transforms, shared by
elements of the same type, with an element-dependent sparse matrix.
This structure reduces the storage and memory access required by the
exact element-block inverse, while sum-factorisation reduces the cost
of applying the transforms. Standard, partial and generalised
sum-factorisation allow the approach to be used on hexahedral, prismatic,
tetrahedral and pyramidal elements within the same framework.
Broadly speaking, we provide a practical guide to constructing 
an efficient matrix-free preconditioner for general DG applications, 
with potential applications beyond the SIPG Poisson solver.

The benchmark results from the SIPG Poisson solver show that efficient preconditioning requires a
balance between convergence and the cost of each iteration.
Point Jacobi is inexpensive but insufficient for strongly coupled modal
bases, whereas block Jacobi improves convergence at a substantial
computational cost. The approximate elemental inverse provides a useful
trade-off between these choices. Although it generally requires more
iterations than the exact block inverse, its lower application cost
can give a shorter solver time. The greatest performance gains are
obtained on high-order hexahedral and prismatic meshes, both as a
standalone preconditioner and within multigrid. The parameter studies
also show the importance of choosing a penalty that ensures stability
without unnecessarily worsening the conditioning.

The gains on tetrahedral and pyramidal meshes are smaller and depend
more strongly on the configuration. For these elements, 
the least-squares fitting is particularly useful
for recovering the approximation quality. However, reducing the local
fitting residual does not guarantee faster global convergence or positive definiteness,
and the resulting approximate inverse is not guaranteed to outperform the exact elemental inverse across all workloads.
Significantly higher setup cost can be another limitation,
although it can be amortised over repeated solves.

The comparison with the continuous Galerkin solver shows that SIPG with
the approximate inverse smoother can achieve similar or shorter solver times on
medium-to-high-order regular hexahedral, prismatic and pyramidal meshes.
Nevertheless, SIPG remains substantially more difficult to solve, particularly
on twisted tetrahedral and pyramidal meshes. 
Its weak enforcement of continuity introduces penalty-dependent
coupling, and the discontinuous space contains more unknowns.
Continuous Galerkin benefits from a much smaller system, better conditioning
and an effective low-energy block preconditioner.
It remains the faster option for tetrahedral and low-order workloads,
and for most of the distorted meshes considered here.

These findings support the motivation for the hybrid multigrid strategy
of \citet{fehn2020hybridMG}, which uses a continuous Galerkin space for
coarse-grid correction while retaining the original discontinuous
problem on the finest level. Transferring to the continuous space
before further coarsening takes advantage of the more efficient
continuous Galerkin discretisation. The current work, 
developing an efficient matrix-free smoother for DG, 
remains valuable in the hybrid framework, as the finest-level smoothing
usually accounts for the majority of the computational cost.

\appendix
\section{Iterative procedure to solve the least-squares problem}
\label{app:fit}

\subsection{Diagonal preconditioning and normal equations}
Eq.~\eqref{eq:fitting-problem} is solved independently on each element.
We therefore omit the subscript $e$ in this discussion for brevity.
Let $\mathbf{t}_i$ denote column $i$ of $\mathbf{T}$. The independent
unknowns are the entries $V_{ij}$ with $(i,j)\in\mathcal{S}$ and
$i\leq j$. The diagonal entries of their least-squares normal matrix are
\begin{equation}
 m_{ij}=
 \begin{cases}
  \|\mathbf{t}_i\|_2^4, & i=j,\\[2pt]
  2\bigl(\|\mathbf{t}_i\|_2^2\|\mathbf{t}_j\|_2^2
    +(\mathbf{t}_i^\T\mathbf{t}_j)^2\bigr), & i<j.
 \end{cases}
 \label{eq:fit-jacobi}
\end{equation}
These values are the squared Frobenius norms of
$\mathbf{t}_i\mathbf{t}_i^\T$ and
$\mathbf{t}_i\mathbf{t}_j^\T+\mathbf{t}_j\mathbf{t}_i^\T$,
respectively. We apply diagonal preconditioning by introducing
scaled unknowns $x_{ij}=\sqrt{m_{ij}}V_{ij}$. Only pairs with $m_{ij}>0$
are active; zero-scale pairs are omitted. The symmetric matrix
$\mathbf{V}(\mathbf{x})$ has entries $V_{ij}=x_{ij}/\sqrt{m_{ij}}$ on
the active pattern and zero elsewhere.

Define the linear map
$\mathcal{F}(\mathbf{x})=\mathbf{T}\mathbf{V}(\mathbf{x})\mathbf{T}^\T$,
whose adjoint operator $\mathcal{F}^*(\cdot)$ is
\begin{equation}
 [\mathcal{F}^*(\mathbf{Y})]_{ij}
 =
 \begin{cases}
  \displaystyle
  \frac{1}{\sqrt{m_{ii}}}\mathbf{t}_i^\T\mathbf{Y}\mathbf{t}_i,
  & i=j, \\
  \displaystyle
  \frac{2}{\sqrt{m_{ij}}}\mathbf{t}_i^\T\mathbf{Y}\mathbf{t}_j,
  & i<j.
 \end{cases}
 \label{eq:fit-adjoint}
\end{equation}
where $\mathbf{Y}\in\mathbb{R}^{N_P\times N_P}$ is an arbitrary symmetric
matrix. 
The least-squares objective becomes $\tfrac12\|\mathbf{I}-\mathcal{F}(\mathbf{x})\|_F^2$,
with normal equations 
\begin{equation}
 \mathcal{F}^*\mathcal{F}(\mathbf{x})
 =\mathcal{F}^*(\mathbf{I}).
 \label{eq:fit-scaled-normal}
\end{equation}
Both operators $\mathcal{F}^*$ and $\mathcal{F}$ are evaluated
through matrix products and restriction to the sparse pattern, without
assembling the least-squares design matrix or its normal matrix.

\subsection{Conjugate-gradient iteration}

\begin{algorithm}[H]
 \caption{Conjugate-gradient iteration for the scaled normal equations}
 \label{alg:fit}
 \begin{algorithmic}[1]
  \Require Projected entries $V_{ij}^{\rm proj}$ from Eq.~\eqref{eq:projected}
           and active scales $m_{ij}>0$. Operators $\mathcal{F}$ and $\mathcal{F}^*$,
           relative tolerance $0<\varepsilon_g<1$, and iteration limit $K_{\max}\geq1$
  \Ensure Approximate least-squares solution $\mathbf{x}$
  \State $x_{ij}^{(0)}\gets\sqrt{m_{ij}}V_{ij}^{\rm proj}$ for each active pair
  \State $\mathbf{R}_0\gets\mathbf{I}-\mathcal{F}(\mathbf{x}_0)$,
         $\mathbf{g}_0\gets\mathcal{F}^*(\mathbf{R}_0)$
  \State $\rho_0\gets\|\mathbf{g}_0\|_2^2$,
         $\mathbf{p}_0\gets\mathbf{g}_0$
  \If{$\rho_0=0$}
   \State \Return $\mathbf{x}_0$
  \EndIf
  \For{$k=0,\ldots,K_{\max}-1$}
   \State $\mathbf{Y}_k\gets\mathcal{F}(\mathbf{p}_k)$,
          $\eta_k\gets\|\mathbf{Y}_k\|_F^2$
   \If{$\eta_k\leq0$}
    \State \Return $\mathbf{x}_k$
   \EndIf
   \State $\alpha_k\gets\rho_k/\eta_k$
   \State $\mathbf{x}_{k+1}\gets\mathbf{x}_k+\alpha_k\mathbf{p}_k$
   \State $\mathbf{R}_{k+1}\gets\mathbf{R}_k-\alpha_k\mathbf{Y}_k$
   \State $\mathbf{g}_{k+1}\gets\mathcal{F}^*(\mathbf{R}_{k+1})$
   \State $\rho_{k+1}\gets\|\mathbf{g}_{k+1}\|_2^2$
   \If{$\rho_{k+1}\leq\varepsilon_g^2\rho_0$ or $k+1=K_{\max}$}
    \State \Return $\mathbf{x}_{k+1}$
   \EndIf
   \State $\beta_k\gets\rho_{k+1}/\rho_k$
   \State $\mathbf{p}_{k+1}\gets\mathbf{g}_{k+1}+\beta_k\mathbf{p}_k$
  \EndFor
  \State \Return $\mathbf{x}_{K_{\max}}$
 \end{algorithmic}
\end{algorithm}

We solve Eq.~\eqref{eq:fit-scaled-normal} by the conjugate-gradient
least-squares method (CGLS) \cite{bjorck1996numerical}, as summarised in
Algorithm~\ref{alg:fit}. This applies conjugate gradient to the scaled
normal equations using only $\mathcal{F}$ and $\mathcal{F}^*$.
The matrix residual is $\mathbf{R}_k=\mathbf{I}-\mathcal{F}(\mathbf{x}_k)$,
and $\mathbf{g}_k=\mathcal{F}^*(\mathbf{R}_k)$ is the residual of the
scaled normal equations. The diagonal preconditioning is already included
in these operators through the scaled unknowns $\mathbf{x}$.

\subsection{Positive-definiteness correction}
\label{app:fit-spd}
The least-squares iteration preserves symmetry but does not ensure
positive definiteness. After the conjugate-gradient iterations, we check the eigenvalues of
$\widetilde{\mathbf{E}}=\mathbf{T}\mathbf{V}\mathbf{T}^\T$, which is
orthogonally similar to $\mathbf{E}_e$ in Section~\ref{sec:quality}, and require
\begin{equation}
 \lambda_{\min}(\widetilde{\mathbf{E}})
 > \epsilon ~,
 \label{eq:fit-spd-check}
\end{equation}
where $\epsilon$ is the chosen tolerance. If this condition is not satisfied, we apply a positive shift $\delta$ to the scaled active diagonal. Let
$\mathbf{D}$ be a diagonal matrix whose entry
$D_{ii}=\|\mathbf{t}_i\|_2^{-1}$ for non-zero columns $i$ and $D_{ii}=0$
otherwise. The corrected matrices are
\begin{equation}
 \begin{aligned}
  \mathbf{V}_*&=\mathbf{V}+\delta\mathbf{D}^2,\\
  \widetilde{\mathbf{E}}_*&=\widetilde{\mathbf{E}}
   +\delta(\mathbf{T}\mathbf{D})(\mathbf{T}\mathbf{D})^\T.
 \end{aligned}
 \label{eq:fit-diagonal-correction}
\end{equation}
We start from $\delta=10^{-6}$ and double it after each check
until a sufficiently large $\delta$ makes $\widetilde{\mathbf{E}}_*$ satisfy Eq.~\eqref{eq:fit-spd-check}.
This empirical correction
can increase the fitting residual and does not minimise the
least-squares objective exactly.

\clearpage
\section{Influence of the penalty scale factor on convergence and accuracy}

\begin{figure}[htbp]
    \centering
    \begin{subfigure}[t]{0.21\textwidth}
        \centering
        \includegraphics[width=\linewidth]{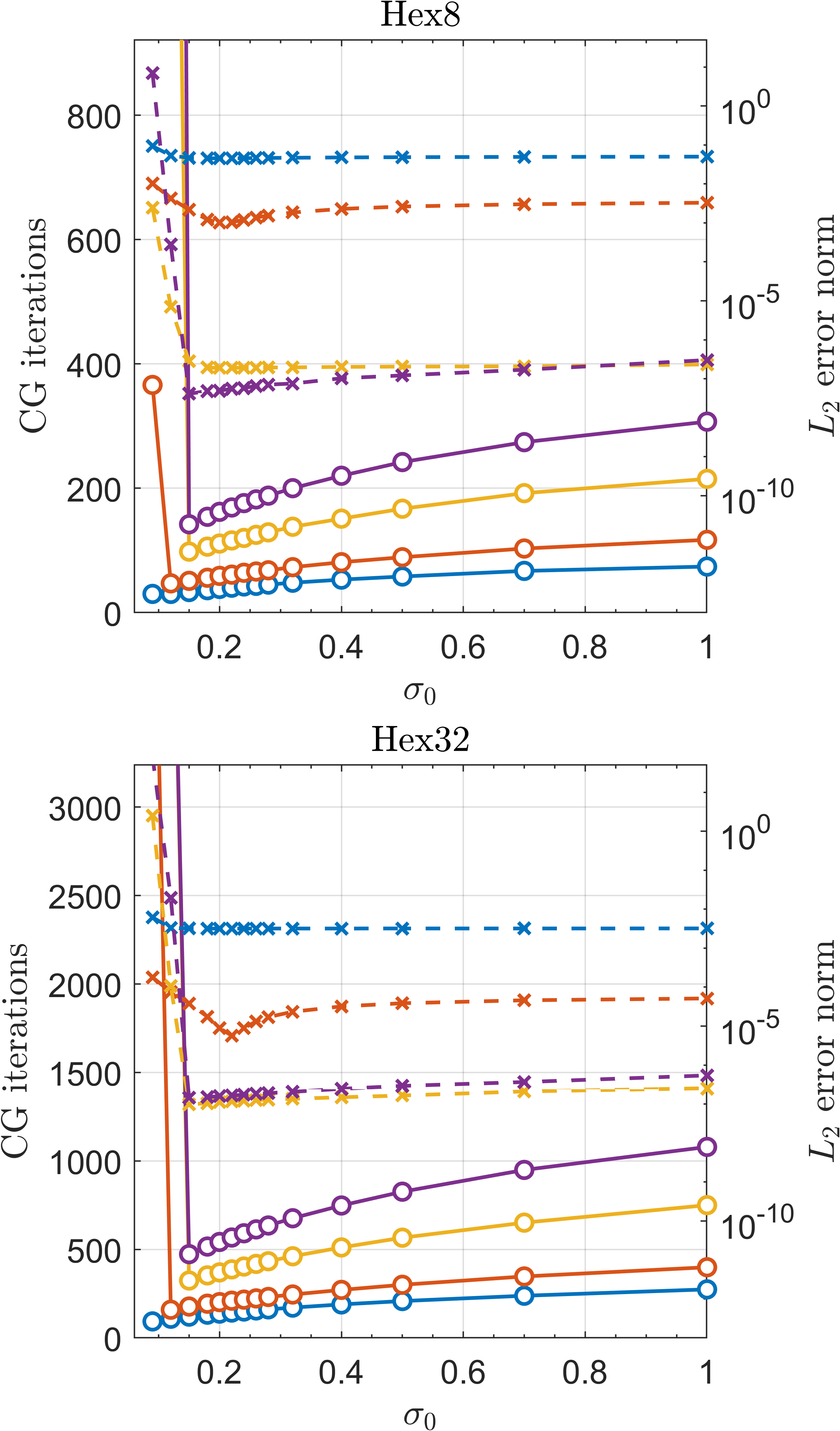}
        \caption{}
        \label{fig:penalty-hex}
    \end{subfigure}
    \begin{subfigure}[t]{0.21\textwidth}
        \centering
        \includegraphics[width=\linewidth]{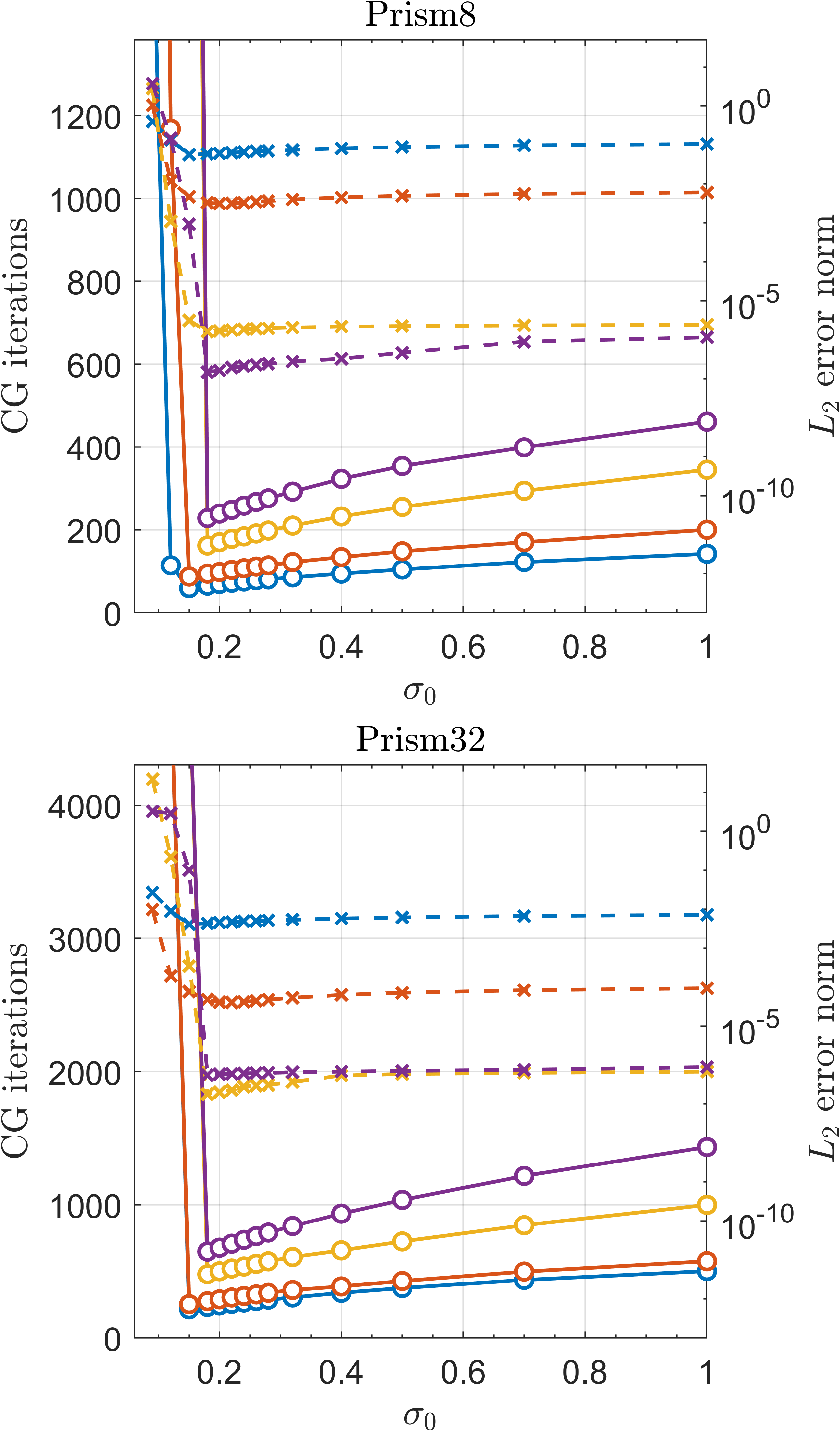}
        \caption{}
        \label{fig:penalty-prism}
    \end{subfigure}
    \begin{subfigure}[t]{0.21\textwidth}
        \centering
        \includegraphics[width=\linewidth]{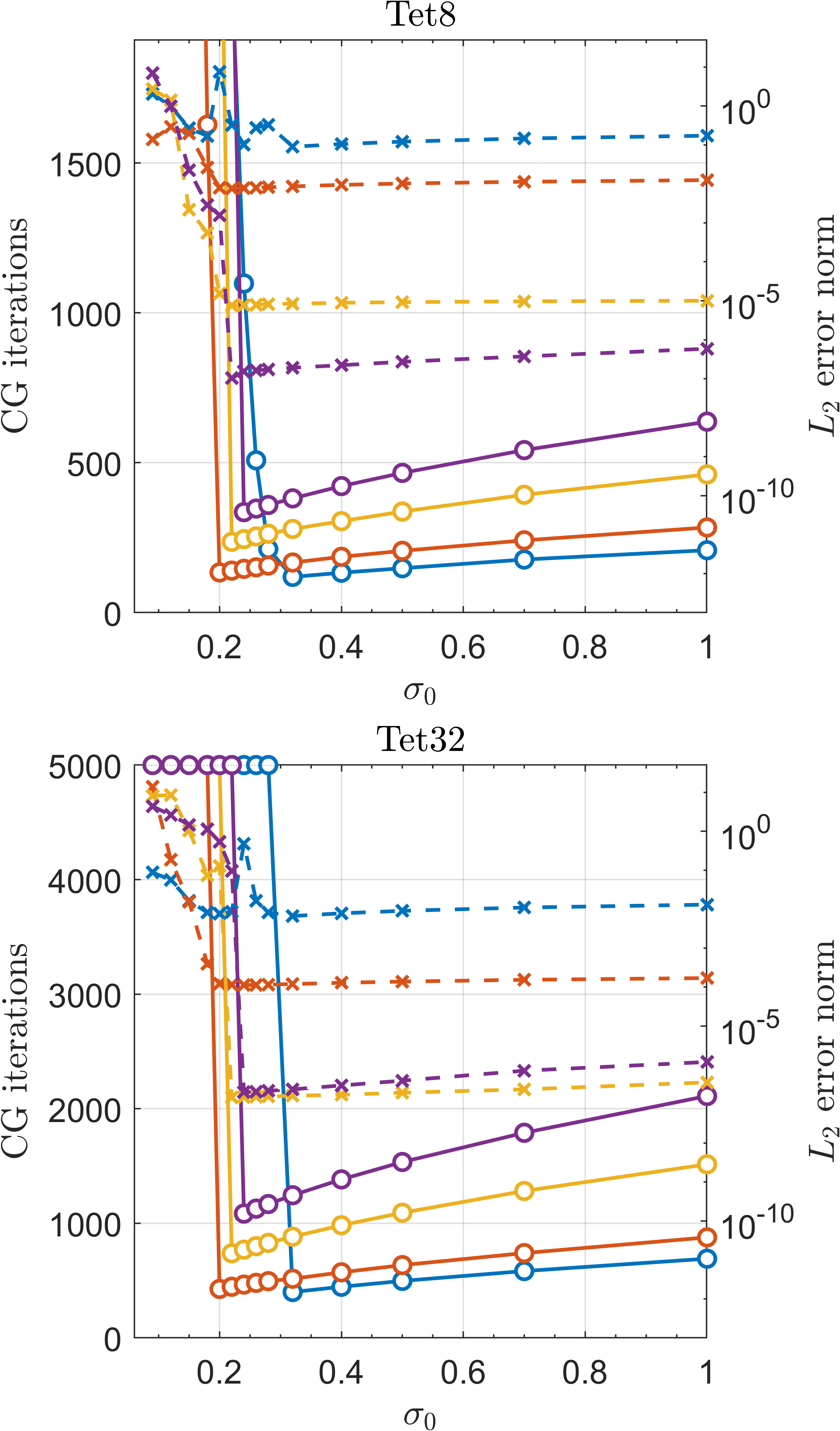}
        \caption{}
        \label{fig:penalty-tet}
    \end{subfigure}
    \begin{subfigure}[t]{0.31\textwidth}
        \centering
        \includegraphics[width=\linewidth]{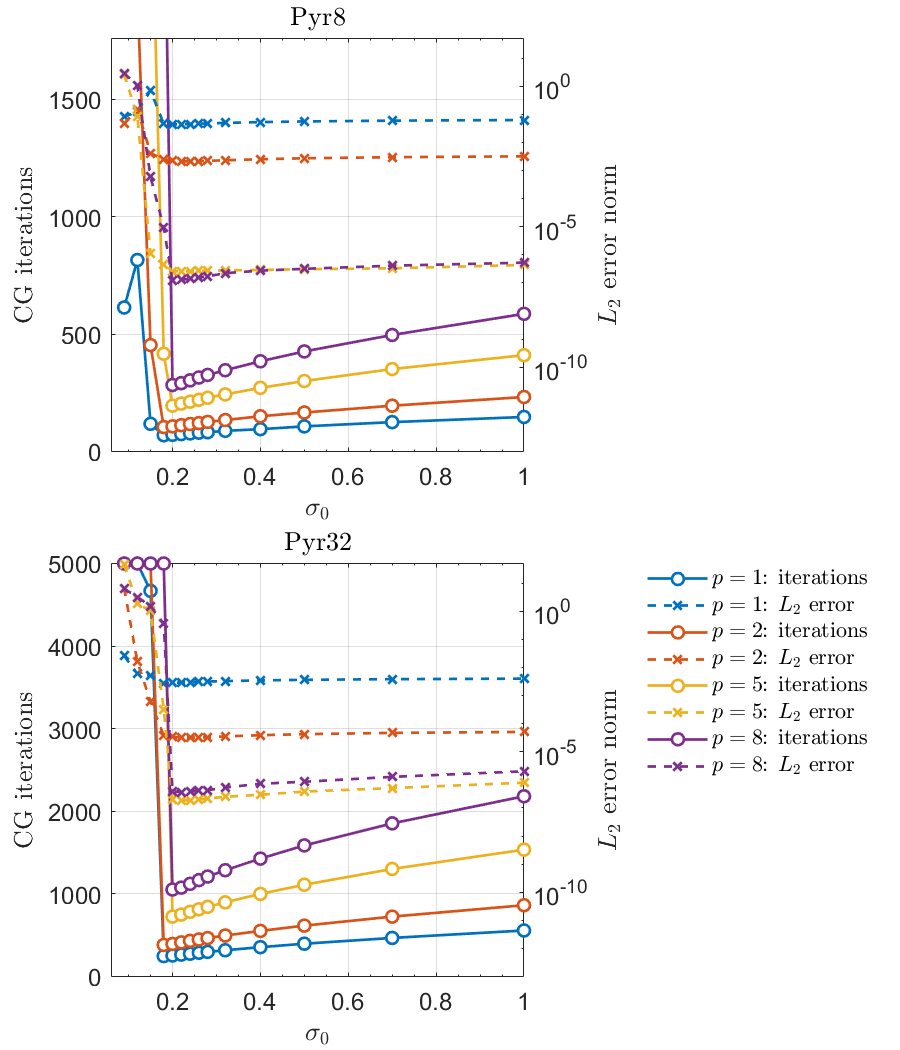}
        \caption{}
        \label{fig:penalty-pyr}
    \end{subfigure}
    \caption{Conjugate gradient iterations and $L_2$ error versus penalty scale factor $\sigma_0$.
        (a) hexahedral mesh; (b) prismatic mesh; (c) tetrahedral mesh; (d) pyramidal mesh.}
    \label{fig:penalty-study-regular}
\end{figure}

\begin{figure}[htbp]
    \centering
    \begin{subfigure}[t]{0.21\textwidth}
        \centering
        \includegraphics[width=\linewidth]{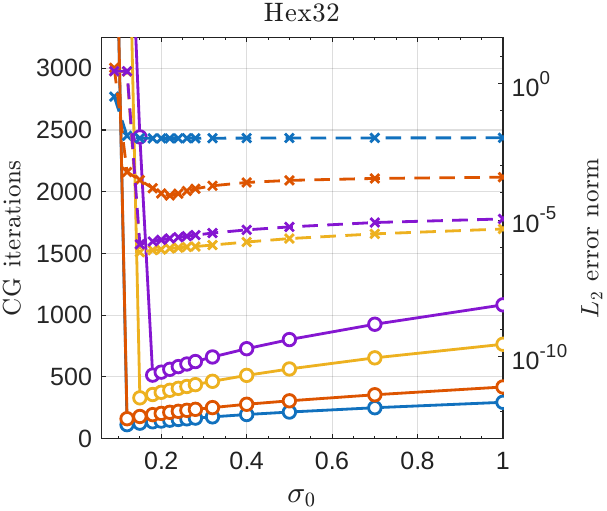}
        \caption{}
        \label{fig:penalty-hex-twisted}
    \end{subfigure}
    \begin{subfigure}[t]{0.21\textwidth}
        \centering
        \includegraphics[width=\linewidth]{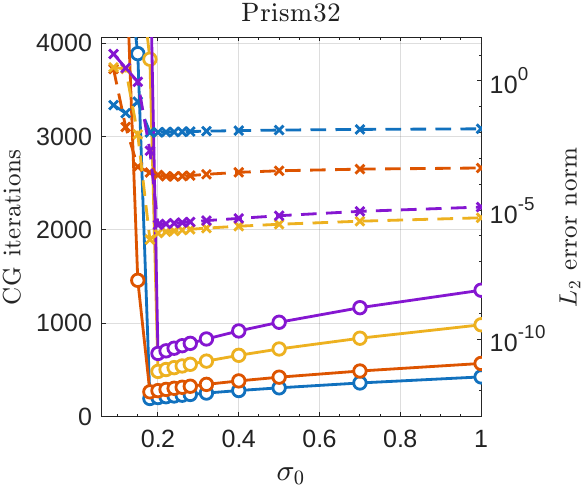}
        \caption{}
        \label{fig:penalty-prism-twisted}
    \end{subfigure}
    \begin{subfigure}[t]{0.21\textwidth}
        \centering
        \includegraphics[width=\linewidth]{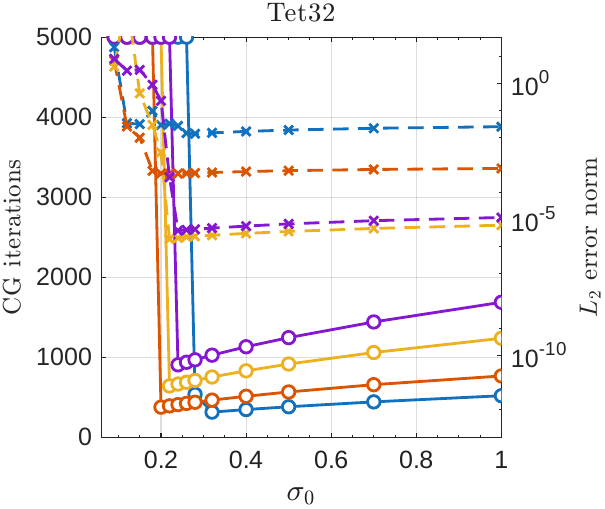}
        \caption{}
        \label{fig:penalty-tet-twisted}
    \end{subfigure}
    \begin{subfigure}[t]{0.31\textwidth}
        \centering
        \includegraphics[width=\linewidth]{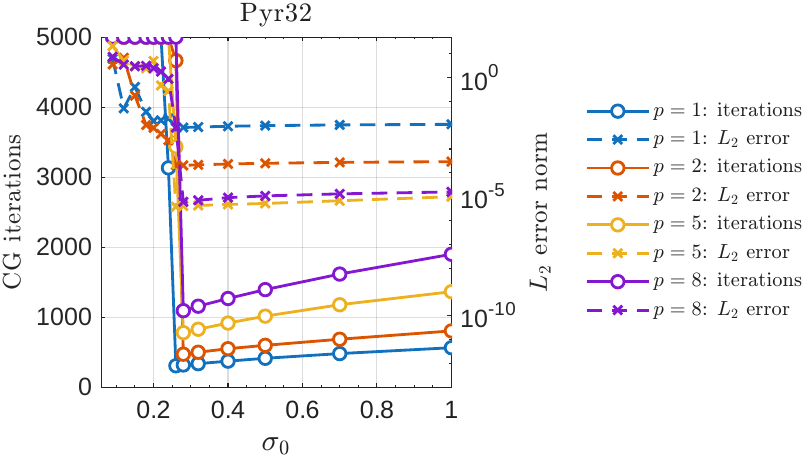}
        \caption{}
        \label{fig:penalty-pyr-twisted}
    \end{subfigure}
    \caption{Conjugate gradient iterations and $L_2$ error versus penalty scale factor $\sigma_0$.
       (a) twisted hexahedral mesh;  (b) twisted prismatic mesh; (c) twisted tetrahedral mesh; (d) twisted pyramidal mesh.}
    \label{fig:penalty-study-twisted}
\end{figure}

\clearpage
\section{Complete preconditioner benchmark results}
\label{app:results}

\begin{figure}[htb]
    \centering
    \includegraphics[width=0.8\textwidth]{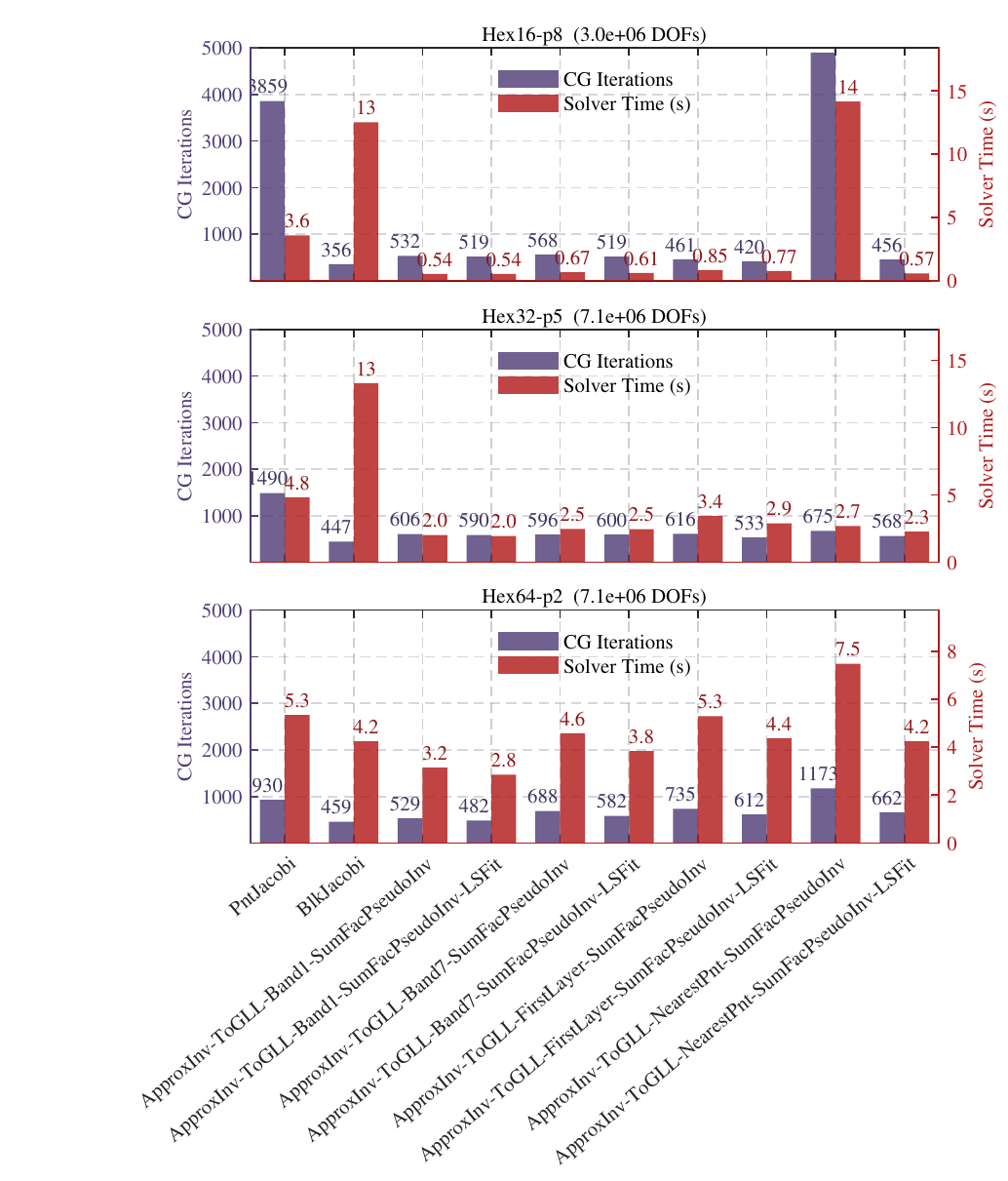}
\caption{The conjugate gradient iterations and solver time required for the SIPG Poisson solver
    to achieve a tolerance of $10^{-9}$ on the \textbf{regular hexahedral} meshes 
    with different preconditioner configurations. The projection space in the approximate inverse
    is fixed to the \textbf{tensorial GLL} grid.}
\label{fig:sipg-preconditioners-hex}
\end{figure}

\begin{figure}[htb]
    \centering
    \includegraphics[width=0.9\textwidth]{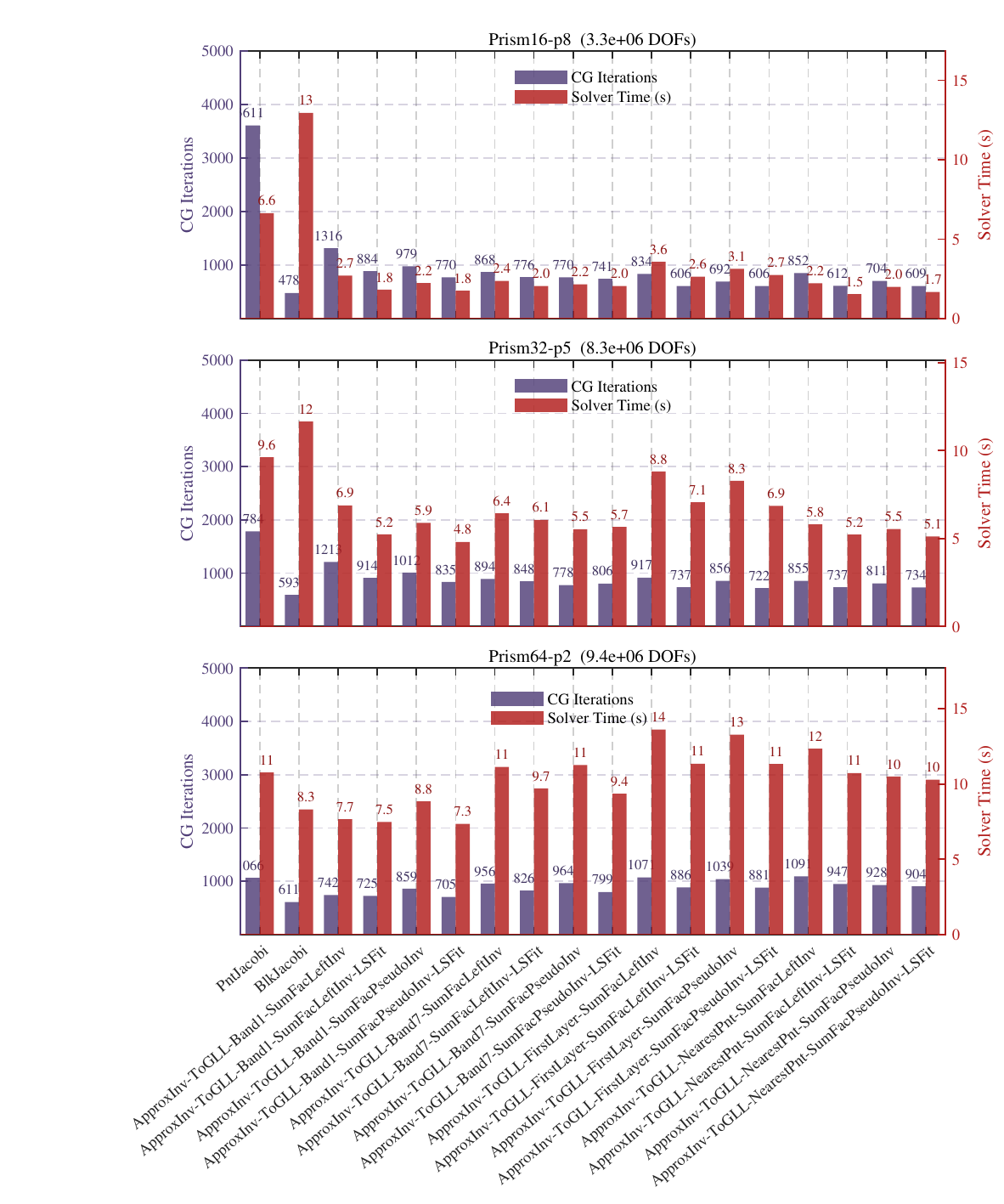}
\caption{The conjugate gradient iterations and solver time required for the SIPG Poisson solver
    to achieve a tolerance of $10^{-9}$ on the \textbf{regular prism} meshes 
    with different preconditioner configurations. The projection space in the approximate inverse
    is fixed to the \textbf{tensorial GLL} grid.}
\label{fig:sipg-preconditioners-prism}
\end{figure}

\begin{figure}[htb]
    \centering
    \includegraphics[width=0.9\textwidth]{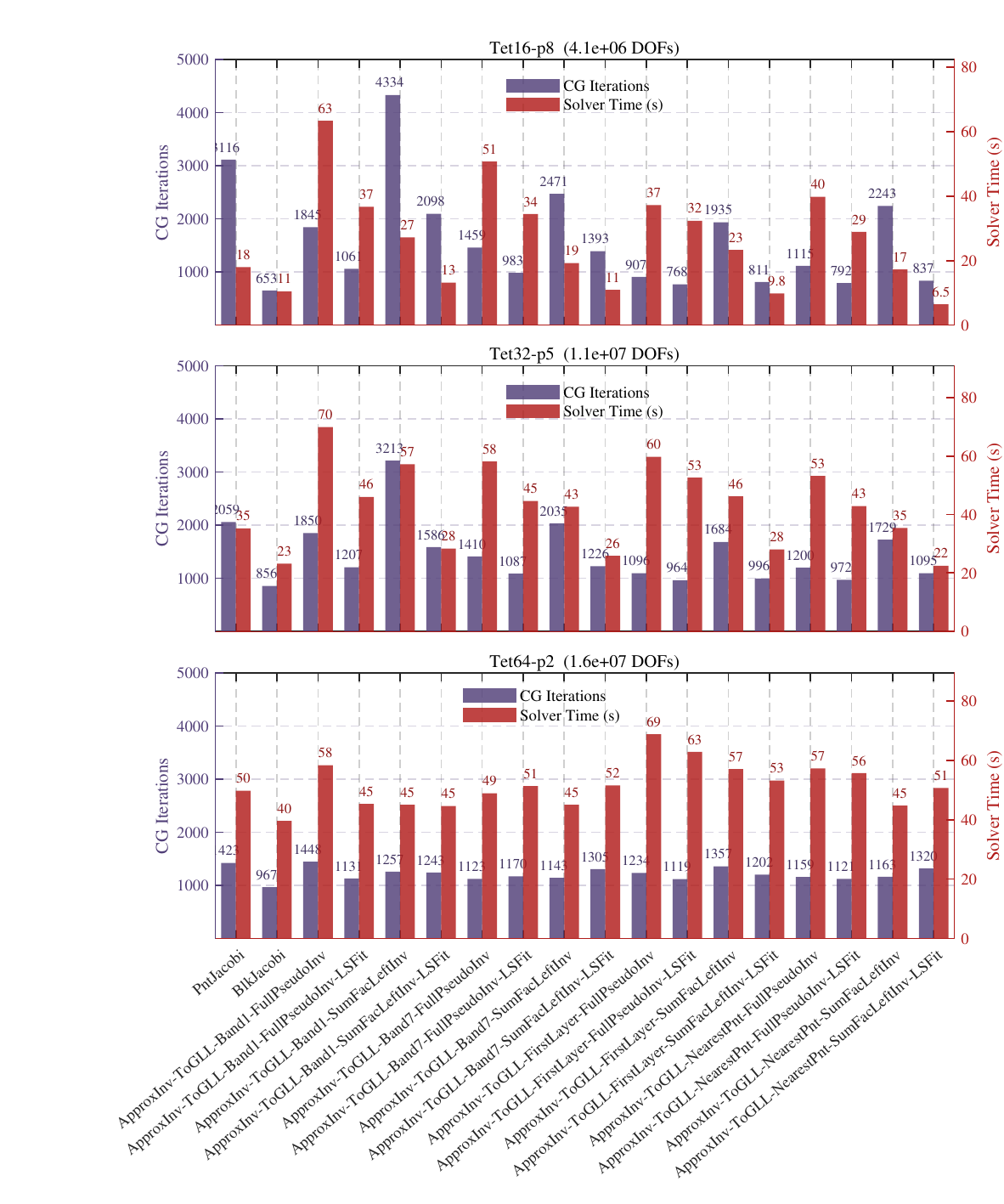}
\caption{The conjugate gradient iterations and solver time required for the SIPG Poisson solver
    to achieve a tolerance of $10^{-9}$ on the \textbf{regular tetrahedral} meshes 
    with different preconditioner configurations. The projection space in the approximate inverse
    is fixed to the \textbf{tensorial GLL} grid.}
\label{fig:sipg-preconditioners-tet}
\end{figure}

\begin{figure}[htb]
    \centering
    \includegraphics[width=0.9\textwidth]{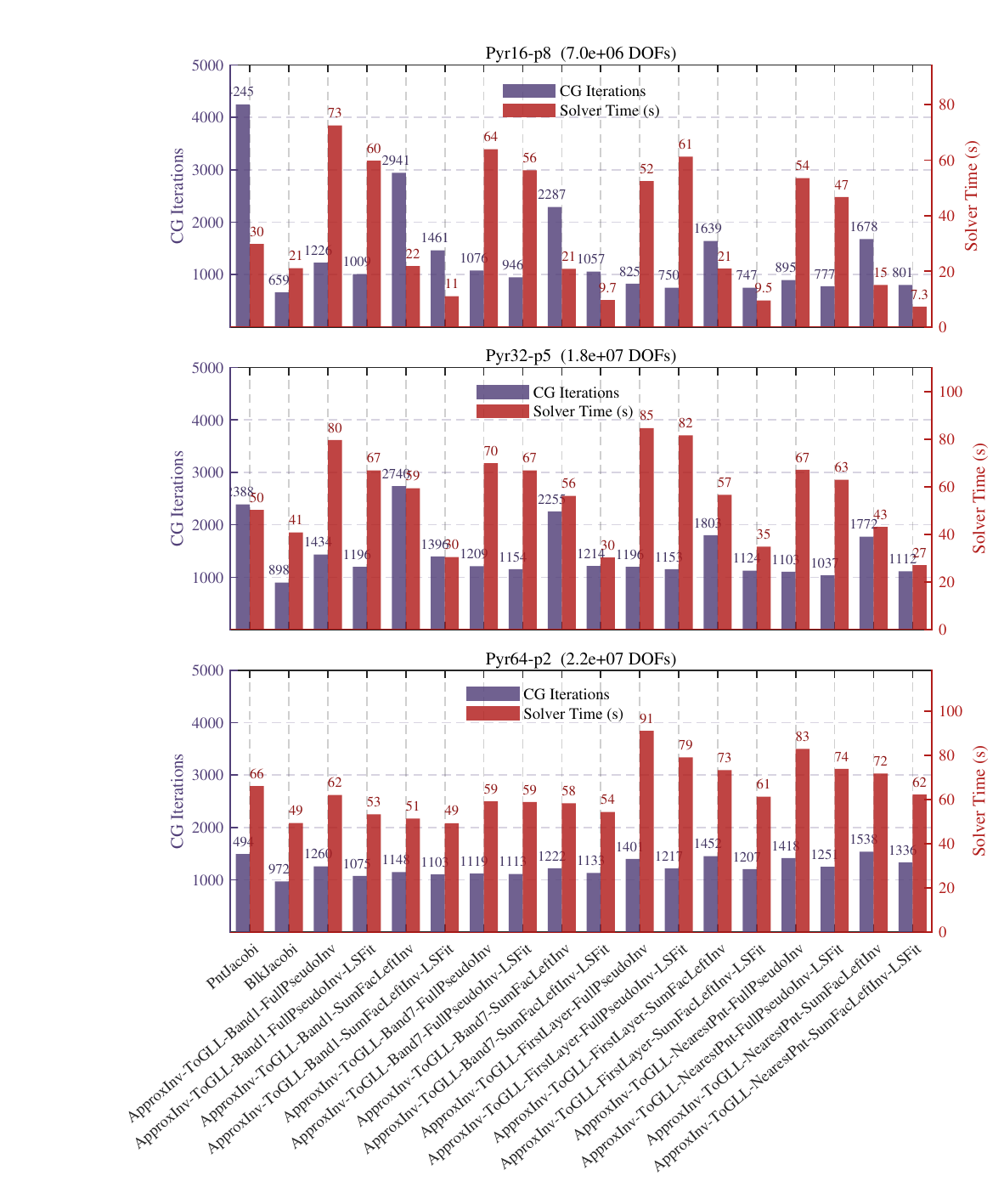}
\caption{The conjugate gradient iterations and solver time required for the SIPG Poisson solver
    to achieve a tolerance of $10^{-9}$ on the \textbf{regular pyramid} meshes 
    with different preconditioner configurations. The projection space in the approximate inverse
    is fixed to the \textbf{tensorial GLL} grid.}
\label{fig:sipg-preconditioners-pyr}
\end{figure}

\begin{figure}[htb]
\centering
    \begin{subfigure}[t]{0.48\textwidth}
        \centering
        \includegraphics[width=\textwidth]{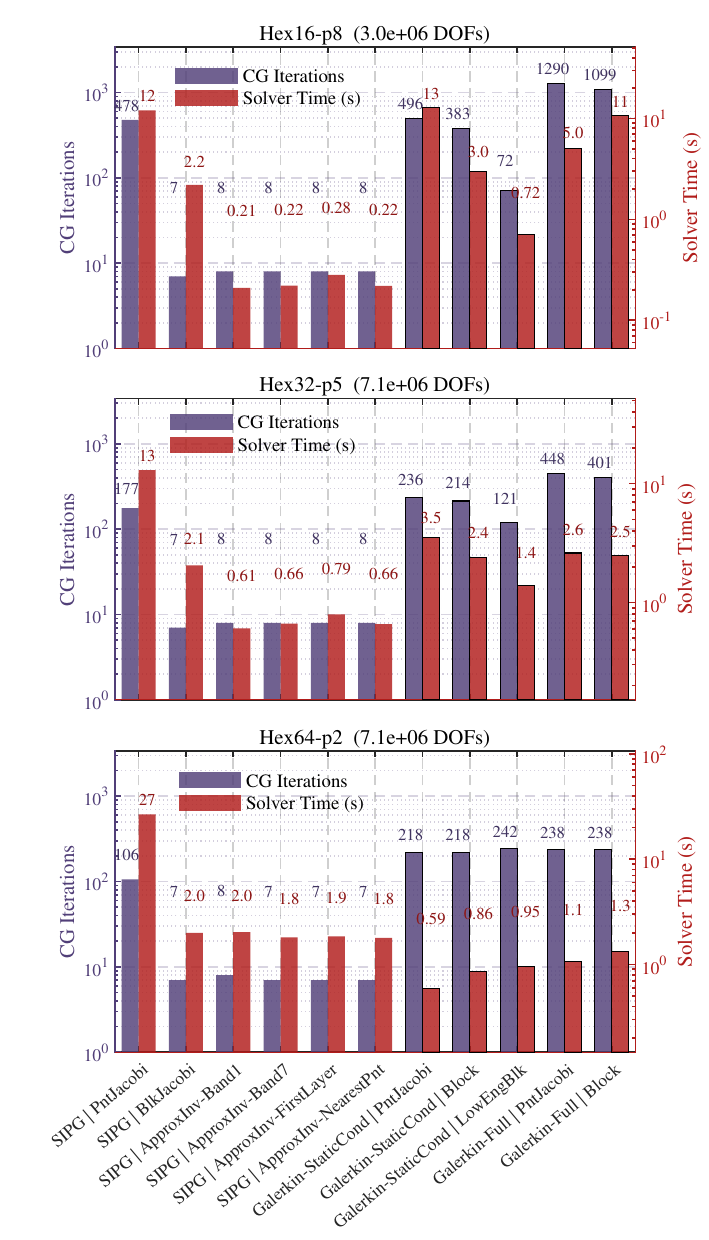}
        \caption{}
    \end{subfigure}
    \begin{subfigure}[t]{0.48\textwidth}
        \centering
        \includegraphics[width=\textwidth]{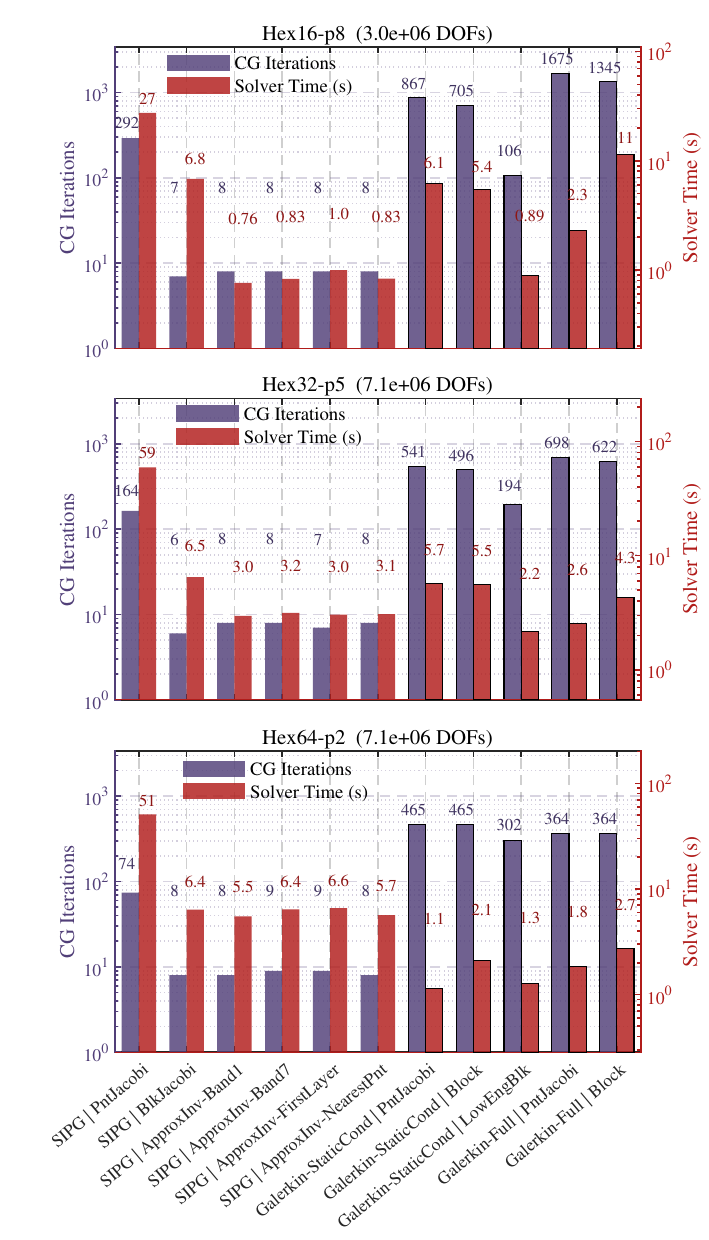}
        \caption{}
    \end{subfigure}
\caption{Comparison of conjugate gradient iteration counts and solver times for the SIPG
    and continuous Galerkin Poisson solvers on the \textbf{(a) regular hexahedral} and \textbf{(b) twisted hexahedral} meshes.
    The SIPG solver employs a $p$-MG preconditioner, and the labels describe different preconditioners inside the Chebyshev smoother. 
    The approximate inverse employs the fastest configurations marked in Table~\ref{tab:different-sparse-pattern}.}
\label{fig:sipg-continuous-hex}
\end{figure}

\begin{figure}[htb]
\centering
    \begin{subfigure}[t]{0.48\textwidth}
        \centering
        \includegraphics[width=\textwidth]{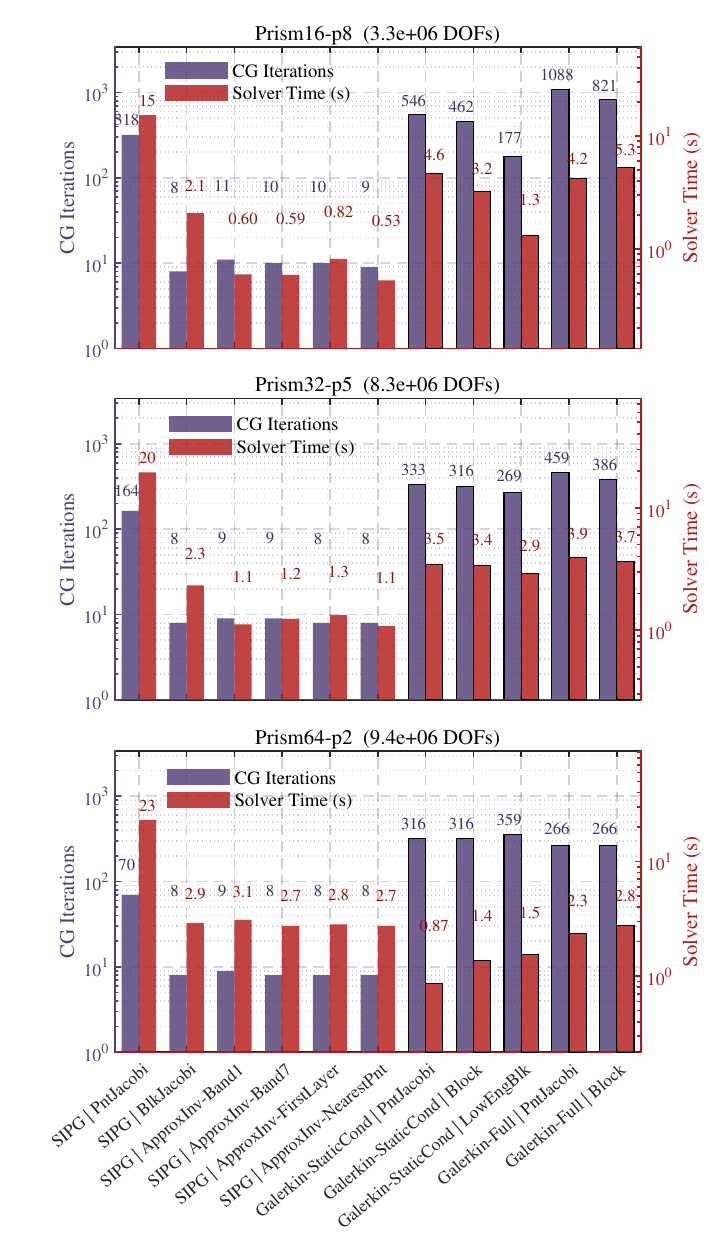}
        \caption{}
    \end{subfigure}
    \begin{subfigure}[t]{0.48\textwidth}
        \centering
        \includegraphics[width=\textwidth]{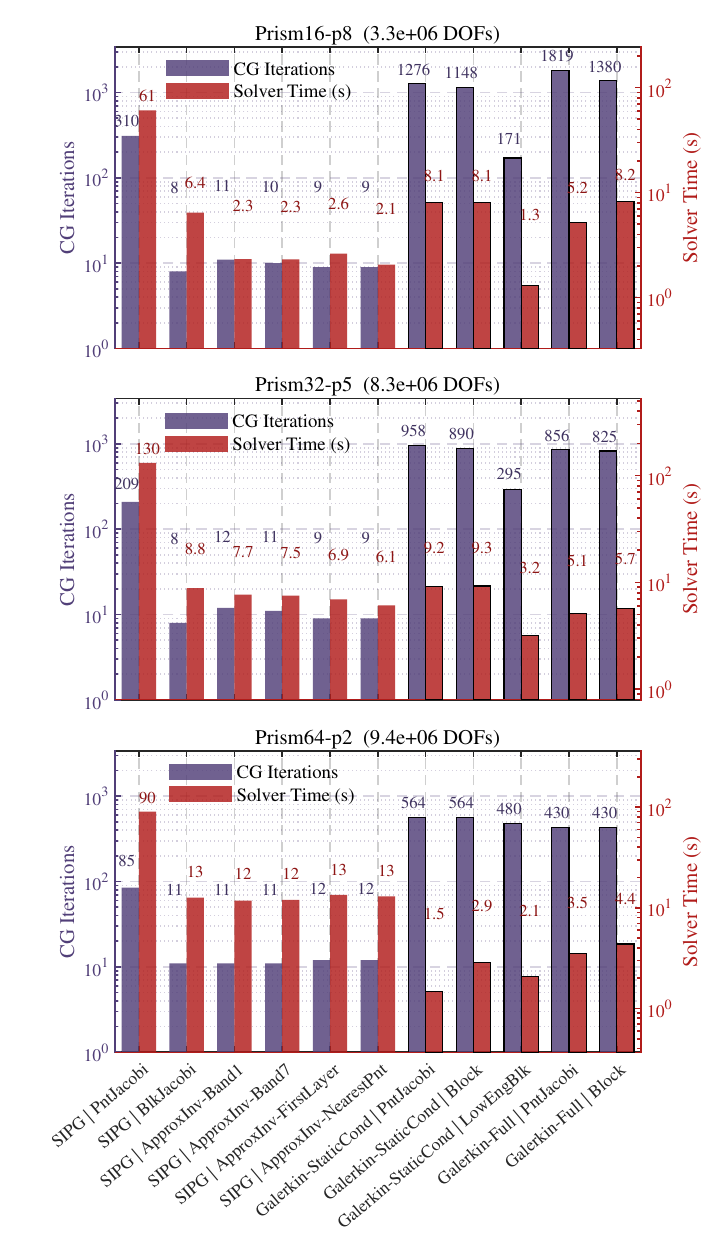}
        \caption{}
    \end{subfigure}
\caption{Comparison of conjugate gradient iteration counts and solver times for the SIPG
    and continuous Galerkin Poisson solvers on the \textbf{(a) regular prismatic} and \textbf{(b) twisted prismatic} meshes.
    The SIPG solver employs a $p$-MG preconditioner, and the labels describe different preconditioners inside the Chebyshev smoother. The approximate inverse employs the fastest configurations given in Table~\ref{tab:different-sparse-pattern}.}
\label{fig:sipg-continuous-prism}
\end{figure}

\begin{figure}[htb]
\centering
    \begin{subfigure}[t]{0.48\textwidth}
        \centering
        \includegraphics[width=\textwidth]{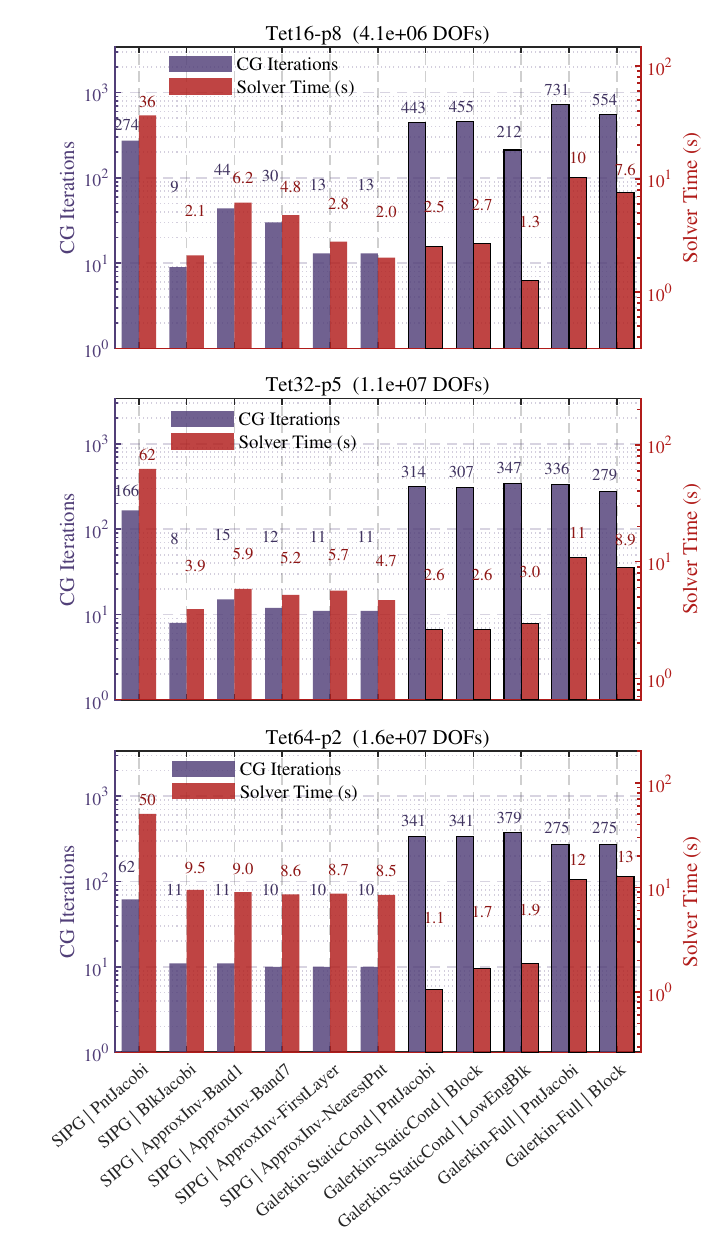}
        \caption{}
    \end{subfigure}
    \begin{subfigure}[t]{0.48\textwidth}
        \centering
        \includegraphics[width=\textwidth]{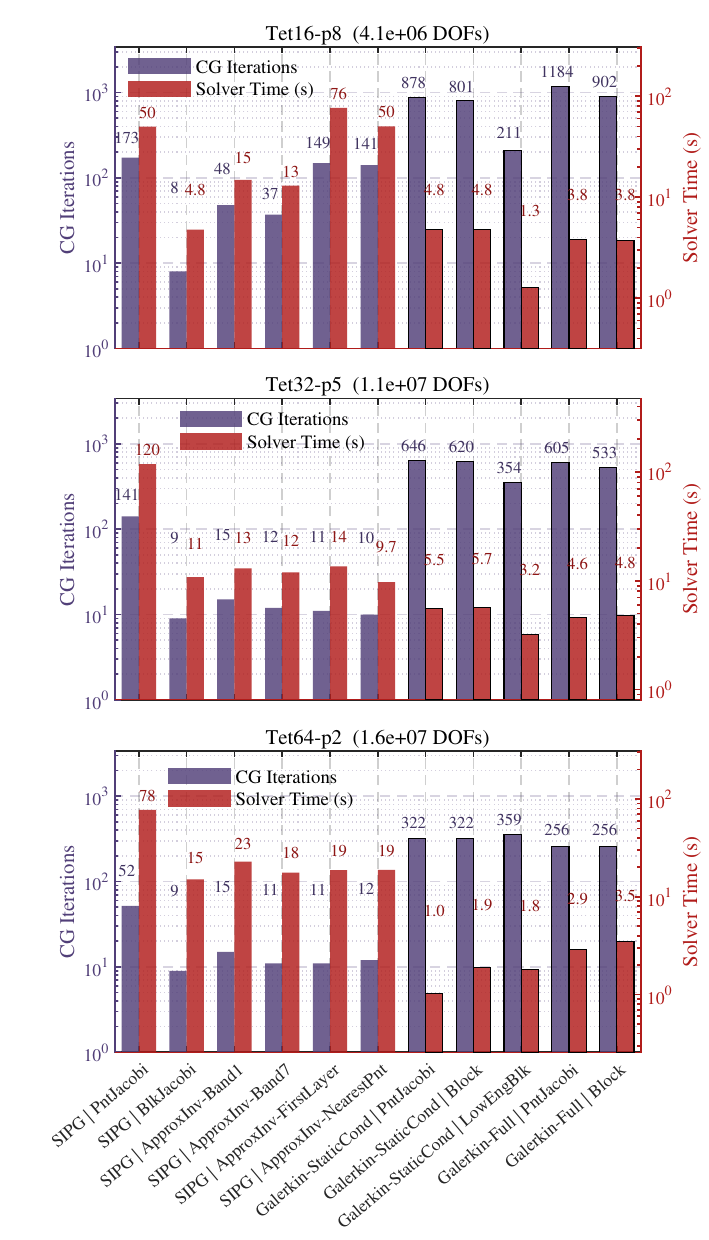}
        \caption{}
    \end{subfigure}
\caption{Comparison of conjugate gradient iteration counts and solver times for the SIPG
    and continuous Galerkin Poisson solvers on the \textbf{(a) regular tetrahedral} meshes and \textbf{(b) twisted tetrahedral} meshes.
    The SIPG solver employs a $p$-MG preconditioner, and the labels describe different preconditioners inside the Chebyshev smoother. The approximate inverse employs the fastest configurations given in Table~\ref{tab:different-sparse-pattern}.}
\label{fig:sipg-continuous-tet}
\end{figure}

\begin{figure}[htb]
\centering
    \begin{subfigure}[t]{0.48\textwidth}
        \centering
        \includegraphics[width=\textwidth]{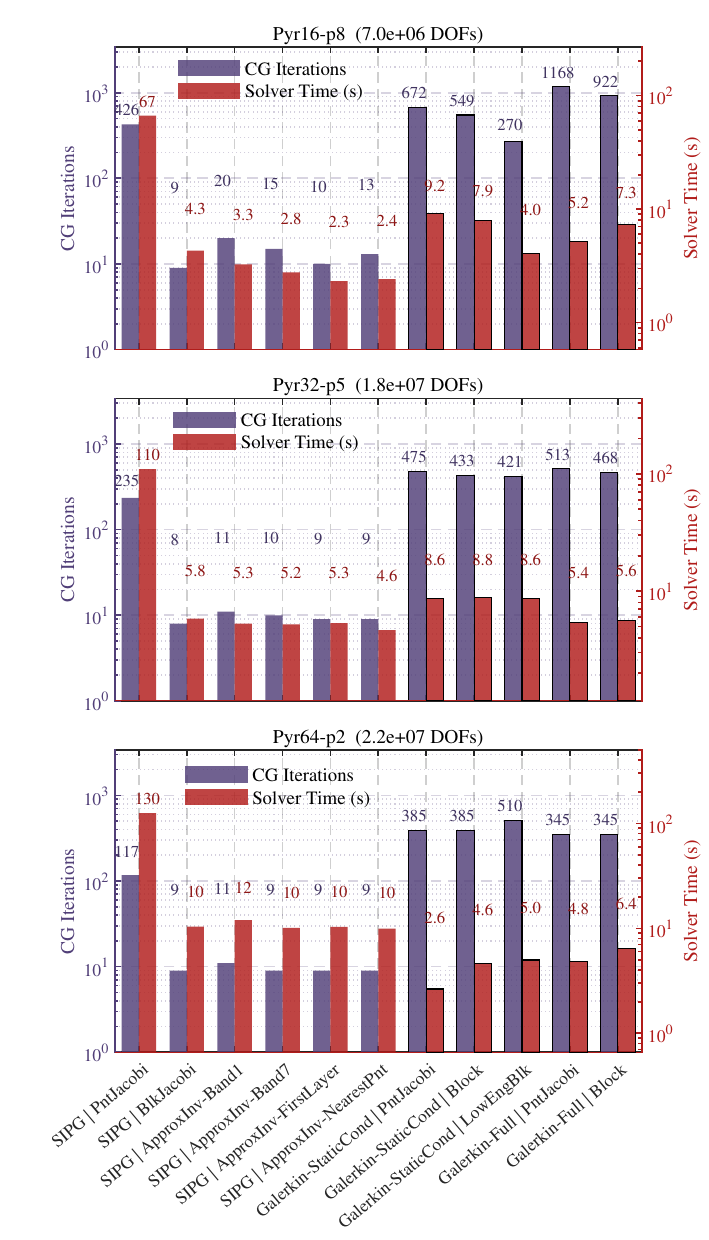}
        \caption{}
    \end{subfigure}
    \begin{subfigure}[t]{0.48\textwidth}
        \centering
        \includegraphics[width=\textwidth]{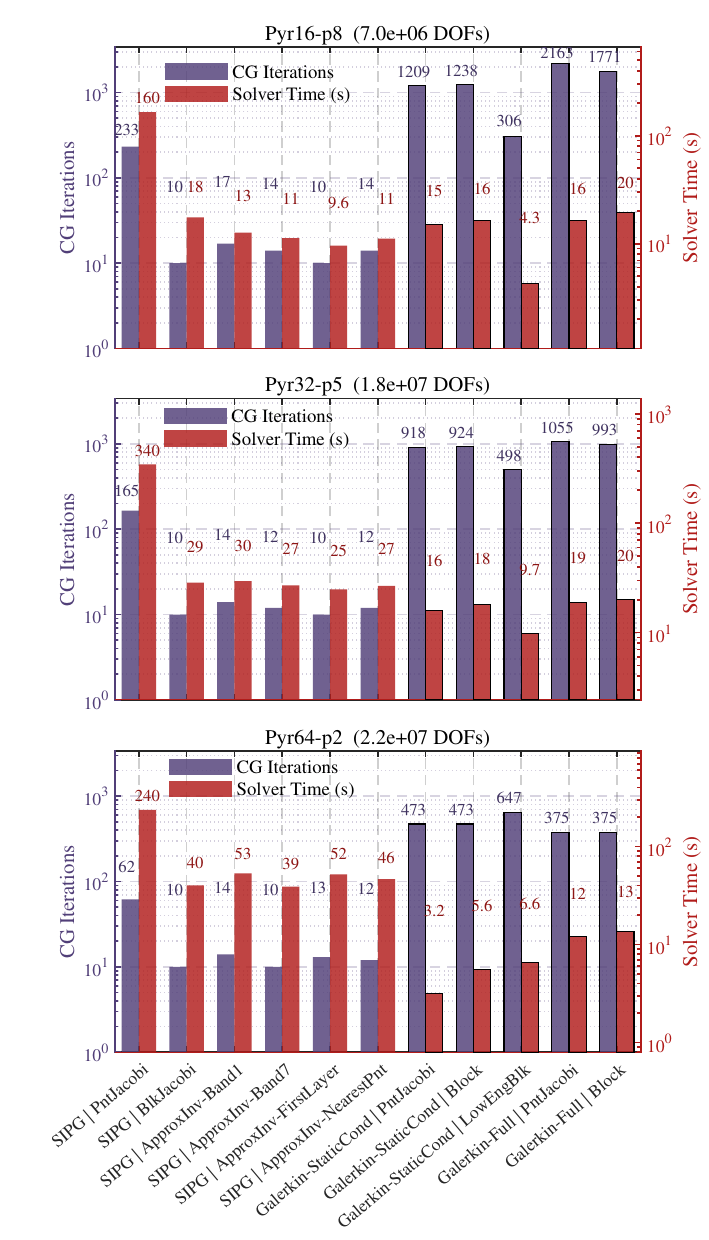}
        \caption{}
    \end{subfigure}
\caption{Comparison of conjugate gradient iteration counts and solver times for the SIPG
    and continuous Galerkin Poisson solvers on the \textbf{(a) regular pyramidal meshes} and \textbf{(b) twisted pyramidal meshes}.
    The SIPG solver employs a $p$-MG preconditioner, and the labels describe different preconditioners inside the Chebyshev smoother. The approximate inverse employs the fastest configurations given in Table~\ref{tab:different-sparse-pattern}.}
\label{fig:sipg-continuous-pyr}
\end{figure}

See Section~\ref{sec:approx-inverse} for the definition of each part of the case label in 
Figs.~\ref{fig:sipg-preconditioners-hex}--\ref{fig:sipg-continuous-pyr}. See Section~\ref{sec:results}
for related discussion.

\clearpage  

\bibliographystyle{unsrtnat}
\bibliography{references}

\end{document}